\documentclass[a4paper,12pt]{article}
\usepackage{epsfig}
\usepackage{amsmath}
\usepackage{amssymb}
\usepackage{amsthm}
\usepackage{bbold}
\usepackage{algorithm2e}
\usepackage{color}
\usepackage{algorithmic}

\usepackage{pgfplots}

\title{CP-TT: using TT-SVD to greedily construct a Canonical Polyadic tensor approximation}
\date{\today}
\author{Virginie Ehrlacher, Maria Fuente Ruiz, Damiano Lombardi}

\begin{document}
\maketitle


\begin{abstract}
  In the present work, a method is proposed in order to compute a Canonical Polyadic (CP) approximation of a given tensor. 
  It is based on a greedy method and an adaptation of the TT-SVD method. The proposed approach can be straightforwardly extended to compute 
  rank-$k$ updates in a stable way. Some numerical experiments are proposed, in which the proposed method is compared to ALS and ASVD methods and performs particularly well for high-order tensors. 
\end{abstract}


\section{Introduction}
\label{sec:intro}
Machine learning and data mining algorithms are becoming increasingly important in analyzing large volume, multi-relational and multi-modal datasets, which are often conveniently represented as multiway arrays or tensors, (\cite{cichocki2016tensor}, \cite{kolda2018tensor}, \cite{kolda2008scalable}). The main challenge in dealing with such data is the so called \textit{curse of dimensionality}, that refers to the need of using a number of degrees of freedom exponentially increasing with the dimension (the reader is refered to \cite{OseTyrt09}). \\

A tensor is said to be in a full format when it is represented as an original multidimensional array. However, distributed storage and processing of high-order tensors in their full format is unfeasible due to the curse of dimensionality. \\

This problem can be alleviated through various distributed and compressed tensor network formats, achieved by low-rank tensor network approximations, as described for instance in \cite{Kolda2009TensorDA}. Some results on low rank decomposition and its uniqueness could  be found in \cite{Os16}, \cite{Domanov2020OnUA}, \cite{Otta14}. It is important to note that, except for very special data structures, a tensor cannot be compressed without incurring in some compression error, since a low-rank tensor representation is only an approximation of the original tensor.

The concept of compression of multidimensional large-scale data by tensor network decompositions can be intuitively explained as follows. Consider the approximation of a $d$-variate function $F(x) =
F(x_1, x_2, . . . , x_d )$ by a finite sum of products of individual functions, each depending on only one or a very few variables. In the simplest scenario, the function $F(x)$ can be (approximately) represented in the following separable form
$$F(x_1,x_2,...,x_d) \approx u^{(1)}(x_1)u^{(2)}(x_2) \cdots u^{(d)}(x_d)$$
In practice, when a $d$-variate function $F(x)$ is discretized into a $d$th-
order array, or a tensor, the approximation above then corresponds to
the representation by rank-1 tensors, also called elementary tensors. 
Denoting by $\mathcal{N}_n$, $n = 1, 2, \cdots, d,$ the size of the discretization grid associated to the $n^{th}$ variate, and by $\mathcal{N} = max_n(\mathcal{N}_n)$, 
the memory requirement to store such a full discretized tensor is equal to $\prod_{n=1}^{d} \mathcal{N}_n \leq \mathcal{N}^{d}$ and grows exponentially with the number of variates $d$. On the other hand, the separable representation of a function is completely defined by its factors, $f^{(n)}(x_n)$, $(n = 1,2,...,d)$, and requires only $\sum_{n=1}^{d} \mathcal{N}_n \ll \mathcal{N}^{d}$ storage units. \\

The separation of variables principle can be achieved trough different tensor formats. In the present work we will make use of two of them: the Canonical Polyadic and the Tensor Train decompositions.

The \textbf{Canonical Polyadic (CP)} decomposition (see \cite{Kolda2009TensorDA}):

Let $r\in\mathbb{N}^*$ be the CP rank, see \cite{DomDeLath14}. The approximation of $F$ in CP format reads:
$$F(x_1,x_2,...,x_d) \approx \sum_{i=1}^r u_i^{(1)}(x_1)u_i^{(2)}(x_2) \cdots u_i^{(d)}(x_d)$$

The \textbf{Tensor Train (TT)} decomposition (see \cite{Oseledets2011TensorTrainD}):
Let $r_1,\ldots r_{d-1}\in\mathbb{N}^*$ be the TT ranks. The approximation reads:
$$F(x_1,x_2,...,x_d) \approx \sum^{r_1}_{i_1=1}... \sum^{r_{d-1}}_{i_{d-1}=1}  u_{i_1}^{(1)}(x_1)u_{i_1, i_2}^{(2)}(x_2)  u_{i_2, i_3}^{(3)}(x_3) \cdots u_{i_{d-2}, i_{d-1}}^{(d-1)}(x_{d-1})  u_{i_{d-1} }^{(d)}(x_d) $$

The main advantages of the CP decomposition are its intuitive expression and its storage scaling. Indeed, instead of the original cost of $\mathcal{O}(\mathcal{N}^d)$, the number of entries 
to store a CP representation reduces to $\mathcal{O}(d\mathcal{N}r)$, which scales linearly in the tensor order $d$ and size $\mathcal{N}$. However, the tensor approximation in CP format may be ill-posed \cite{deSilvLim08} and leads to numerical instabilities. The compression of tensors using the CP is usually computed by means of the ALS (alternating least squares) method, which sometimes has some performance issues \cite{Beylkin05algorithmsfor} especially for high-order tensors. Some modifications have been proposed in order to make the compression more efficient, as in \cite{8081290}, \cite{els}. 

The Tensor Train format is probably one of the most used tensor formats in realistic applications (\cite{bigoni2016,wang2019distributed,rakhuba2016}), due to a good trade off between optimality and numerical stability. 
The TT format combines two advantages to take into consideration: on the one hand, it is stable from an algorithmic point of view; on the other, it is computationally affordable provided that the TT ranks of the tensors stay reasonably small. \\
The number of entries in the TT format is $\mathcal{O}(rd\mathcal{N}+r^3(d-2))$. Even if the number of entries could be larger than in CP,  the main advantage of the TT format is its ability to provide stable quasi-optimal rank reduction, obtained, for instance, by truncated singular value decompositions. 
 
In the literature, hybrid formats combining CP with other methods have been proposed in \cite{Os20}, and described in \cite{grasedyck2013literature}. Also, CP has been 
combined with TT in \cite{kour2020efficient}, where we can see from a different perspective, that both methods combined have potential improvements.\\

The main contribution of the present work is a method that constructs a CP tensor approximation by exploiting the TT-SVD algorithm. We 
will refer to it as CP-TT. The proposed method is based on a greedy algorithm in which we reduce the approximation to a sequence of rank-1 (or rank-$k$ for any $k\in \mathbb{N}^*$) approximations. 
The method, although relying on the TT-SVD iteration does not require to fix \emph{a priori} the order of the variables. 

The work is structured as follows: in Section \ref{sec:method} the generic greedy strategy to define a 
CP decomposition is presented along with two alternatives to the proposed method, namely ALS and ASVD
(respectively in Sections \ref{subsec:ALS} and \ref{subsec:ASVD}). The formulation of the present approach is
presented in Section \ref{subsec:CP-TT}. The properties of the method, the computational cost and the extension to compute stable rank$-k$ updates 
are presented in Section \ref{sec:math}. Some numerical experiments and results are presented in Section \ref{sec:numExp}.

\section{The method}
\label{sec:method}
Before detailing the method formulation and its properties, we introduce hereafter the notation and the problem setting. Then, different methods are recalled, which were proposed in the literature in order to construct a CP tensor approximation of a given tensor. Among them, in the present work, we focus on ALS and on ASVD. The Section ends with the presentation of the CP-TT method.

\subsection{Notation}
\label{subsec:Notation}
The notation is presented. The method proposed in the present work was motivated by high-dimensional function approximation. We therefore make the choice to deal with continuous tensors. The extension to discrete tensors is straightforward.

Let $d\in\mathbb{N}^*$ be the tensor order. Let $i=1,\ldots,d$, and $d_i\in\mathbb{N}^*$. We introduce $\Omega^{(i)}\in\mathbb{R}^{d_i}$, which are open bounded sets. Let the domain be denoted by: $\Omega = \Omega^{(1)} \times \ldots \times \Omega^{(d)}$: the variables are $x\in\Omega$, $x=(x_1,\ldots, x_d)$.

A real valued tensor $F$ is a function defined as:
\begin{equation}\label{eq:datum}
F : \left\{ \begin{array}{ll}
                                                 \Omega \ \rightarrow \ \mathbb{R} \\
                                                   x \ \mapsto F(x_1,\ldots, x_d).\\
                                                 \end{array}\right.
\end{equation} 
Let $k\in\mathbb{N}^*$, and $1\leq p\leq \infty$. In the present work we restrict to the problem of approximating a function $F\in L^2(\Omega)$. Observe that $L^2(\Omega)= L^2(\Omega^{(1)})\otimes\ldots\otimes L^2(\Omega^{(d)})$.
From now on, we denote the $L^2$ scalar product of two nonzero functions $u \in L^2(\Omega)$ and $v \in L^2(\Omega)$ :
$$
\langle u, v \rangle = \int_{\Omega} u v \ dx, 
$$
and the $L^2$ norm as $\|u\|_{L^2(\Omega)}^2=\langle u, u \rangle$.\\

A pure tensor product is a function defined as follows:
\begin{equation}\label{eq:tensor}
u^{(1)}\otimes\ldots\otimes u^{(d)} : \left\{ \begin{array}{ll}
                                                 \Omega^{(1)}\times\ldots\Omega^{(d)} \ \rightarrow \ \mathbb{R} \\
                                                   (x_1,\ldots, x_d) \ \mapsto u^{(1)}(x_1)\ldots u^{(d)}(x_d).\\
                                                 \end{array}\right.
\end{equation}

Let $n\in\mathbb{N}^*$, $c_i\in \mathbb{R}, \ i=1,\ldots,n$ be a sequence of real numbers, a tensor $T$ in CP format is written as follows:
\begin{equation}
\label{eq:CP}
T = \sum_{i=1}^n c_i u^{(1)}_i\otimes\ldots\otimes u^{(d)}_i. \nonumber
\end{equation} 
The number of terms $n$ is the CP-rank of the tensor. The set of the rank$-n$ tensors in CP format is denoted as:
\begin{equation}
\mathcal{T}^{(n)} = \left \lbrace T = \sum_{i=1}^n c_i u^{(1)}_i\otimes\ldots\otimes u^{(d)}_i, \ c\in\mathbb{R}^n, \ u_i^{(j)}\in L^2(\Omega_j) \ \forall i=1, \ldots, n, j=1,\ldots, d \right\rbrace
\end{equation}

The objective is to construct a tensor approximation of a given function. The problem reads as follows:
\begin{eqnarray}\label{eq:best}
T_*^{(n)} = \mathrm{arg}\inf_{T\in \mathcal{T}^{(n)}} \| F - T \|^2_{L^2(\Omega)}.
\end{eqnarray}
A solution $T_*^{(n)}$ to (\ref{eq:best}) is called a best rank$-n$ approximation of $F$. This problem is known to be ill-posed, in the sense that there may not exist any solutions to (\ref{eq:best}). 
This translates into numerical instabilities for most algorithms aiming at computating CP approxiamtions of tensors~\cite{de2008tensor}. 

One way to approach the solution of the problem is to adopt a greedy strategy, that consists in replacing the best rank$-n$ approximation by a sequence of rank$-1$ approximations, which are well defined.
The general strategy reads as follows:
\begin{eqnarray*}
T^{(0)} = 0, \\
t^{(k+1)}_* \in \mathrm{arg}\inf_{t\in\mathcal{T}^{(1)}} \| F - T^{(k)} - t\|^2_{L^2(\Omega)}, \\
T^{(k+1)} = T^{(k)} + t_*^{(k+1)}. 
\end{eqnarray*} 

\medskip

There are several algorithms in order to compute a rank$-1$ approximation solution to the above problem, and the methods presented hereafter aim at computing a numerical approximation of a solution $t_*^{(k+1)}$ of the above minimization problem.

\subsection{Alternating least square (ALS)}
\label{subsec:ALS}

The Alternating Least Square (ALS) method is based on fixed point method iterations and is one of the most used method to compress tensors in CP format. 
It consists in the following iterative scheme to compute an approximation $t_{\rm ALS}^{(k+1)}$ of $t_*^{(k+1)}$, for a given error threshold $\epsilon >0$. 

\bfseries ALS: \normalfont

\begin{itemize}
 \item \bfseries Initialization: \normalfont 
 Choose randomly $u^{(1,0)} \in  L^2(\Omega^{(1)})$, ...,  $u^{(d,0)} \in  L^2(\Omega^{(d)})$ and set $m\geq 1$. 
 
 \item \bfseries Iteration $m\geq 1$:  \normalfont 
 
\begin{enumerate}
 
 \item For $i = 1, \cdots , d$, compute $u^{(i,m)} \in  L^2(\Omega^{(i)})$ solution to 
 $$
 u^{(i,m)} \in \mathop{{\rm argmin}}_{u^{(i)}\in L^2(\Omega^{(i)})} \| F - T^{(k)} - u^{(1,m)} \otimes u^{(2,m)} \otimes u^{(i-1,m)} \otimes u^{(i)} \otimes u^{(i+1, m-1)} \otimes \cdots \otimes u^{(d, m-1)} \|^2_{L^2(\Omega)}. 
 $$
 
 \medskip
 
 \item Compute $\eta:= \left\| u^{(1,m)} \otimes \cdots \otimes u^{(d, m)} - u^{(1,m-1)} \otimes \cdots \otimes u^{(d, m-1)} \right\|_{L^2(\Omega)}$. 
 
 \medskip
 
 \item If $\eta < \epsilon$, define $t_{\rm ALS}^{(k+1)}:= u^{(1,m)} \otimes \cdots \otimes u^{(d, m)}$ the approximation of the best rank-one approximation of $F - T^{(k)}$. Else, set $m =2$ and iterate again.  
 \end{enumerate}
\end{itemize}

\medskip

The convergence properties of the ALS algorithm have been abundantly studied. We refer the reader for more details to the following series 
of works~\cite{uschmajew2012local,rohwedder2013local,espig2015convergence,wang2015accelerating,oseledets2018alternating}.

\subsection{Alternating Singular Value Decomposition (ASVD)}
\label{subsec:ASVD}

For the sake of comparison with the ALS method and our proposed procedure, we introduce in this section the so-called Alternating Singular Value Decomposition (ASVD) method, which was proposed in~\cite{friedland2013best}.

We denote by $\mathcal I:= \left\{ \{i,j \}, \; 1\leq i < j \leq d\right\}$ be the set of all possible pairs of indices between $1$ and $d$. An ordering of the elements of $\mathcal{I}$ is chosen so that 
$$
\mathcal I= (I_k)_{1\leq k \leq K}
$$
where $K = |\mathcal I|$.

The ASVD method then also consists in an iterative scheme to compute an approximation $t_{\rm ASVD}^{(k+1)}$ of $t_*^{(k+1)}$, for a given error threshold $\epsilon >0$, which reads as follows: 

\bfseries ASVD: \normalfont

\begin{itemize}
 \item \bfseries Initialization: \normalfont 
 Choose randomly $u^{(1,0)} \in  L^2(\Omega^{(1)})$, ...,  $u^{(d,0)} \in  L^2(\Omega^{(d)})$ and set $m\geq 1$. 
 
 \item \bfseries Iteration $m\geq 1$:  \normalfont 
 
\begin{enumerate}
 \item Set $u^{(i,m)} = u^{(i,m-1)}$ for all $1\leq i \leq d$. 
 
 \medskip
 
 \item For $k = 1, \cdots , K$, let $I_k = (i_k, j_k)$ and compute $U^{(k,m)} \in  L^2(\Omega^{(i_k)} \times \Omega^{(j_k)})$ solution to 
 $$
 U^{(k,m)} \in \mathop{{\rm argmin}}_{U^{(k)}\in L^2(\Omega^{(i_k)} \times \Omega^{(j_k)} )} \| F - T^{(k)} - U^{(k)} \otimes \bigotimes_{1\leq l \leq d; l\neq i_k, j_k}u^{(l,m)} \|^2_{L^2(\Omega)}. 
 $$
 Update $u^{(i_k, m)}$ and $u^{(j_k, m)}$ so that $(u^{(i_k, m)},u^{(j_k, m)}) \in L^2(\Omega^{(i_k)}) \times L^2(\Omega^{(j_k)})$ is solution to 
  $$
 (u^{(i_k,m)},u^{(j_k,m)})  \in \mathop{{\rm argmin}}_{ (u^{(i_k)},u^{(j_k)})\in L^2(\Omega^{(i_k)}) \times L^2(\Omega^{(j_k)} )} \| U^{(k,m)} - u^{(i_k)}\otimes u^{(j_k)} \|^2_{L^2(\Omega)}. 
 $$
 
 \medskip
 
 \item Compute $\eta:= \left\| u^{(1,m)} \otimes \cdots \otimes u^{(d, m)} - u^{(1,m-1)} \otimes \cdots \otimes u^{(d, m-1)} \right\|_{L^2(\Omega)}$. 
 
 \medskip
 
 \item If $\eta < \epsilon$, define $t_{\rm ASVD}^{(k+1)}:= u^{(1,m)} \otimes \cdots \otimes u^{(d, m)}$ the approximation of the best rank-one approximation of $F - T^{(k)}$. Else, set $m =2$ and iterate again.  
 \end{enumerate}
 
\end{itemize}

\subsection{CP-TT}
\label{subsec:CP-TT}
In this section, the CP-TT method is introduced. The idea behind this method of tensor approximation is to combine the CP format and the TT-SVD method, in order to benefit from the simplicity of the CP format and the numerical stability of the TT-SVD. 

\begin{figure}[htbp]
\begin{center}
\includegraphics[scale=0.4]{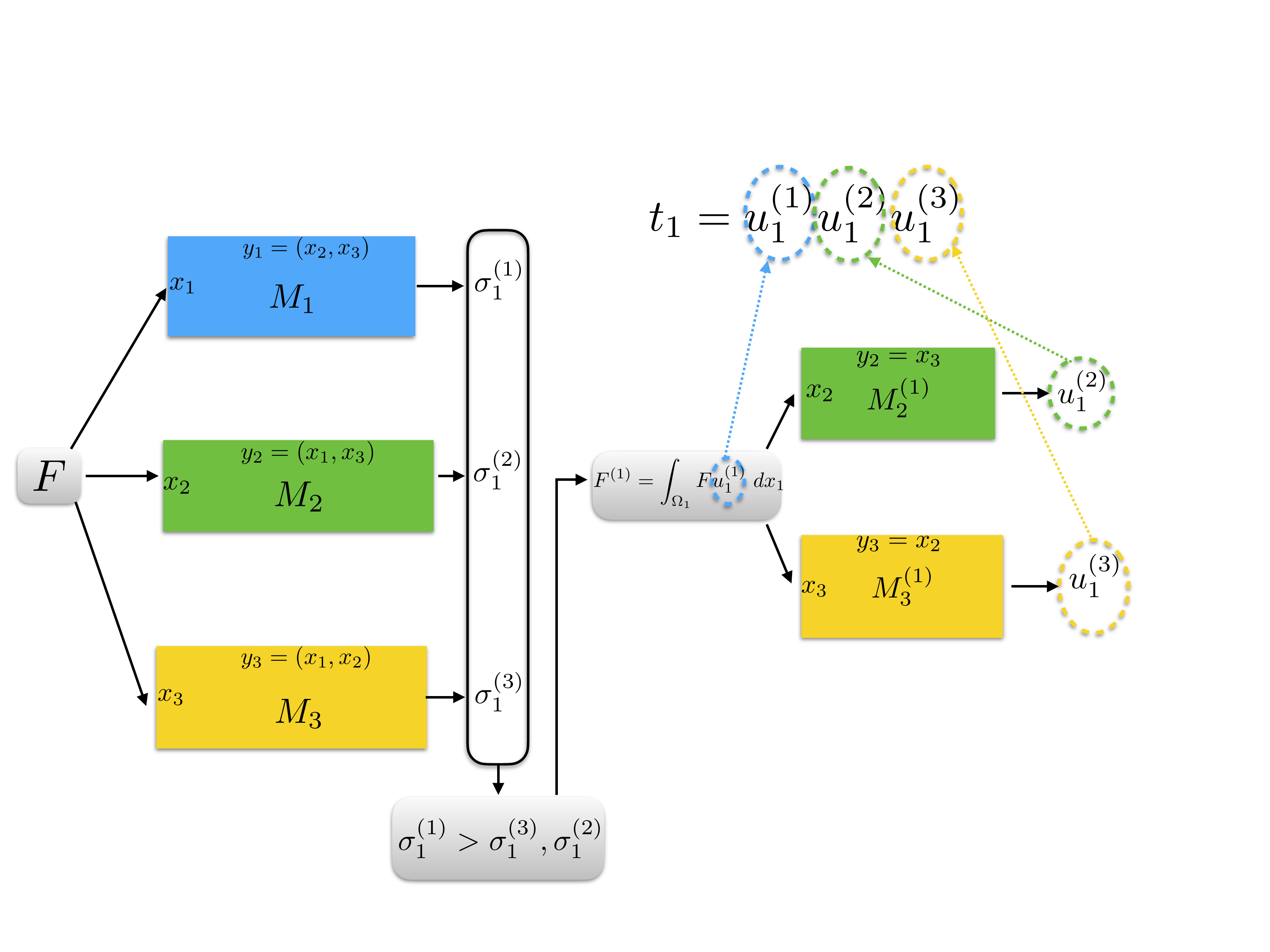}
\caption{Schematic representation of the CP-TT iteration for a $3$rd order tensor $F$.}
\label{cpttscheme}
\end{center}
\end{figure}

The principle of the method (which is sketched in Fig.\ref{cpttscheme} in the particular case of a $3$rd order tensor) relies itself on a greedy algorithm, which is detailed below. Without loss of generality, 
for the sake of simplicity, we present the method in full details in the case where $k=0$, i.e. 
when $T^{(0)} =0$ so that $F - T^{(0)} = F$.

For all $1 \leq i \leq d$, let us define $\widehat{\Omega}_i:=\Omega_1 \times \cdots\times \Omega_{i-1} \times \Omega_{i+1} \times \cdots \times \Omega_d$ and let us introduce $M_i: \Omega_i \times \widehat{\Omega}_i \to \mathbb{R}$ the $i^{th}$ unfolding of $F$ defined by
$$
M_i(x_i,y_i) = F(x_1,\cdots,x_d)
$$
for all $x_i\in \Omega_i$ and all $y_i=(x_1,\cdots,x_{i-1},x_{i+1},\cdots,x_d)\in \widehat{\Omega}_i$.

Let us now consider a POD decomposition of $M_i$:
\begin{equation}
M_i(x_i,y_i) = \sum_{j>0} \sigma_j^{(i)}u_j^{(i)}(x_i) v_j^{(i)}(y_i),
\end{equation}
where $u_j^{(i)}$ is an orthonormal basis of $L^2(\Omega_i)$, $v_j^{(i)}$ is an orthonormal basis of $L^2(\widehat{\Omega}_i)$ and $\sigma_j^{(i)}\in\mathbb{R}^+$ are the singular values, that we assumed to be ranked in decreasing order. 

We define $1\leq i_* \leq d$ as the integer such that:
\begin{equation}
i_* = \mathrm{arg}\max_{1\leq i\leq d} \sigma_1^{(i)}.
\end{equation}

A rank-one approximation of $M_{i_*}$, and therefore of $F$, is then constructed as $m_{i_*} = \sigma_{1}^{(i_*)}u_1^{(i_*)}\otimes v_1^{(i_*)} = u_1^{(i_*)} \otimes w_1^{(i_*)}$, where $w_1^{(i_*)}\in L^2(\widehat{\Omega}_{i_*})$ is defined as: $w_1^{(i_*)}=\sigma_1^{(i_*)}v_1^{(i_*)}$. \\

Without loss of generality (up to changing the labeling of the variables), let us assume that $i_* = 1$. It then holds that: 
$$w_1^{(1)}(y_1) =\int_{\Omega_1}M_1(x_1, y_1)u_1^{(1)}(x_1) dx_1, \ \ \forall y_1 = (x_2,\ldots,x_d) \in \widehat{\Omega}_1$$

We introduce an auxiliary tensor, defined as follows:

$$F^{(1)}(x_2, \ldots, x_d) = w_1^{(1)}(y_1) = \int_{\Omega_1}F(x_1,\ldots,x_d) u_1^{(1)}(x_1) dx_1.$$

In the next step of the method, we look for an approximation of $F^{(1)} $, by proceeding in an analogous way. 

For all $2 \leq i \leq d$, let us define $\widehat{\Omega}_i^{(1)}:=\Omega_2 \times \cdots\times \Omega_{i-1} \times \Omega_{i+1} \times \cdots \times \Omega_d$ and let us introduce $M_i^{(1)}: \Omega_i \times \widehat{\Omega}_i^{(1)} \to \mathbb{R}$ the $i^{th}$ unfolding of $F^{(1)}$ defined by:
$$
M_i^{(1)}(x_i,y_i^{(1)}) = F^{(1)}(x_2,\cdots,x_d)
$$
for all $x_i\in \Omega_i$ and all $y_i^{(1)}=(x_2,\cdots,x_{i-1},x_{i+1},\cdots,x_d)\in \widehat{\Omega}_i^{(1)}$.

A POD decomposition of $M_i^{(1)}$ is computed, of the form:
\begin{equation}
M_i^{(1)}(x_i,y_i^{(1)}) = \sum_{j>0} \sigma_j^{(i,1)}u_j^{(i,1)}(x_i) v_j^{(i,1)}(y_i^{(1)}),
\end{equation}
where $u_j^{(i,1)}$ is an orthonormal basis of $L^2(\Omega_i)$, $v_j^{(i,1)}$ is an orthonormal basis of $L^2(\widehat{\Omega}_i^{(1)})$ and $\sigma_j^{(i,1)}\in\mathbb{R}^+$ are the singular values, that we assumed to be ranked in decreasing order. 

We define $2\leq i_*^{(1)} \leq d$ as the integer such that:
\begin{equation}
i_*^{(1)} = \mathrm{arg}\max_{2\leq i\leq d} \sigma_1^{(i,1)}.
\end{equation}

Analogously to what we have presented above, up to changing the labeling of the variables, we can assume that $i_*^{(1)}=2$. We define:
\begin{equation}
F^{(2)}(x_3,\ldots,x_d) = \int_{\Omega_2} F^{(1)}(x_2,\ldots, x_d) u^{(2,1)}_1(x_2) \ d x_2.
\end{equation}

This process continues until it remains a $2-$nd order tensor and the POD decomposition gives us the final terms. Let us suppose, for sake of simplicity in the notation, that the two last variables are $x_{d-1}, x_{d}$. The truncated rank-one approximation of $F^{(d-2)}$ is given by:
$$F^{(d-2)}(x_{d-1},x_d) \approx \sigma_1^{(d-1, d-2)} u_1^{(d-1, d-2)}(x_{d-1}) v_1^{(d, d-2)}(x_d).$$

The final rank-one approximation of $F$ is obtained as:

$$t^{(1)}_{CPTT} = \sigma_1^{(d-1, d-2)} u_1^{(1)} \otimes u_1^{(2,1)} \otimes \dots \otimes  u_1^{(d-1, d-2)}\otimes v_1^{(d, d-2)}.$$

\paragraph*{Remark:} We have considered, for sake of simplicity in the notation, the order $I_* = \left\lbrace 1,2,\ldots, d \right\rbrace$. We stress here that in the proposed method, the order is not fixed \emph{a priori}. Instead, it is the result of the iteration, and, in particular, of the optimization step on $\sigma_i^{(\cdot)}$. 

The CP-TT method is summarized in the pseudocode \ref{alg:algCPTT}. 

\subsection{Optimization of the CP coefficients}
To obtain a more accurate approximation, we introduce an optimization problem in order to find the best linear combination of the rank-1 tensor obtained by one of the three methods presented above. Let us denote 
by $t^{(1)}, \cdots, t^{(k)}$ the $k$ pure rank-1 tensors obtained after $k$ iterations of each of the three methods mentioned above (ALS, ASVD or CPTT).

For each method, at the $k^{th}$ iteration, we look for coefficients $c:=(c_1,\ldots,c_k) \in\mathbb{R}^k$ such that:
$$T^{(k)}(c)= \sum_{j=1}^{k}c_j t^{(j)}(x_1,...,x_d)$$
minimizes the norm of the residual. The optimization problem reads:
 $$(c_{*,1}, \cdots, c_{*,k}) =:c_{*}  =  \mathrm{arg} \min_{c \in \mathbb{R}^k}\left(\frac{1}{2} \| F-T^{(k)}(c) \|_{L^2(\Omega)}^2\right) $$

We define the functional $J$:
$$
J: \left\{
\begin{array}{ccc}
\mathbb{R}^{k}  & \to & \mathbb{R}^+\\
 c & \mapsto & \|F-T^{(k)}(c)\|^2_{L^2(\Omega)}\\
\end{array} \right.
$$
   
The Euler Lagrange equations read as follows: 
  
  $$
  \forall 1\leq l \leq k, \quad 
  \frac{\partial J}{\partial c_l}(c_*) =\int_{\Omega} (F-\sum_{j=1}^k c_{*,j} t^{(j)})(-t^{(l)})dx_1\ldots dx_d=0$$
This implies that
$$  
 \forall 1\leq l \leq k, \quad   \int_\Omega F t^{(l)} dx_1\ldots dx_d = \sum_{j=1}^k\int_\Omega  c_{*,j} t^{(j)} t^{(l)} dx_1\ldots dx_d,$$
which reduces to the following linear system
$$
A c_* = b
$$
where $A = (A_{lj})_{1\leq l,j \leq k} \in \mathbb{R}^{k\times k}$, $b = (b_l)_{1\leq l \leq k} \in \mathbb{R}^k$ with
$$
\forall 1\leq l,j \leq k, \quad A_{lj}= \int_\Omega t^{(j)} t^{(l)} dx_1\ldots dx_d \mbox{ and } b_l = \int_\Omega F t^{(l)} dx_1\ldots dx_d.
$$

\medskip

For each of the three methods mentioned above (ALS, ASVD, CP-TT), the resulting procedure then reads as follows: 

\bigskip

\bfseries ALGORITHM 1: \normalfont
\medskip
\begin{enumerate}
 \item \bfseries Initialization: \normalfont $k=0$, $T^{(0)} = 0$; 
 \item \bfseries Iteration $k$: \normalfont Compute $t^{(k+1)}$ an approximate solution to the optimization problem
 $$
\mathrm{arg}\inf_{t\in\mathcal{T}^{(1)}} \| F - T^{(k)} - t\|^2_{L^2(\Omega)},
 $$
 using either the ALS, ASVD or CP-TT method. 
 
 Compute $$
 c_*^{(k+1)}:= \mathrm{arg} \min_{c \in \mathbb{R}^{k+1}}\left(\frac{1}{2} \| F-T^{(k+1)}(c) \|_{L^2(\Omega)}^2\right). 
 $$
 Define $T^{(k+1)}:= \sum_{l=1}^{k+1} c_{*,l}^{(k+1)} t^{(l)}$ and set $k:= k+1$.
\end{enumerate}

%
%
%
%
%
%
\begin{algorithm}
\caption{Scheme of the computation algorithm of the CP-TT term}
\begin{algorithmic}
\medskip
 
\STATE{\textbf{The CP-TT method:}}
\medskip
\STATE{We define $F^{(0)}:=F$}
\STATE{$I^{(0)}=\{1,\dots,d\}$}
\STATE{Set $m=1$}
\WHILE{$\# I^{(m)} > 2$}
      \STATE{Compute the unfoldings $(M_i)_{i\in I^{(m-1)}}$ of $F^{(m-1)}$}
      \STATE{Compute their POD approximation, and denote by $(\sigma_j^{i})_{j\geq 1}$ the singular values of $M_i$ for all $i\in I^{(m-1)}$}   
      \STATE{Select the largest singular value: $i_{*,m}=\mathrm{arg}\max_{i\in I^{(m-1)}} \sigma_1^{i}$ }
      \STATE{Define $I^{(m)}:= I^{(m-1)} \setminus \{i_{*,m}\}$}
      \STATE{Select $u^{i_{*,m}}_1$ the first singular function associated to the POD decomposition of $M_{i_{*,m}}$}
      \STATE{Define $F^{(m+1)}:= \int_{\Omega_{i_*}} F^{(m)} u^{i_{*,m}}_1 dx_{i_*}$}
      \STATE{Set $m:= m+1$.}
 \ENDWHILE
\STATE{Let $I^{(d-2)} = \{ i_{*,d-1}, i_{*,d}\}$. Compute the POD decomposition of $F^{(d-2)}$}
\STATE{Let $u^{i_{*,d-1}}_1$ and $u^{i_{*,d}}_1$ be the first singular functions of $F^{(d-2)}$}
\STATE{Define $t^{(1)}_{\rm CPTT} = \bigotimes_{j=1}^d u^{i_{*,j}}_1$.} 
\end{algorithmic}
\label{alg:algCPTT}
\end{algorithm}

\section{Properties of the CP-TT method}
\label{sec:math}
In this section, several properties of the proposed CP-TT method are analyzed. 

\subsection{Orthogonality properties}

Let us analyze the inner iteration used in order to compute a rank-1 update. For sake of simplicity, without loss of generality we consider the sequence 
of $i_*$ corresponding to the largest singular values is $1,2, \ldots, d$. Let us define

$$G^{(1)}(x_1,\ldots, x_d) := F- u_1^{(1)}(x_1)F^{(1)}(x_2, \ldots, x_d).$$
Then, it holds that
$$\|G^{(1)}\|_{L^2(\Omega)}^2 = \sum_{j>1}\left(\sigma_{j}^{(1)}\right)^2. $$

By projecting the residual $F$ onto $u_1^{(1)}$ we get 

\begin{equation}
 \int_{\Omega_1}  u_1^{(1)}(x_1)F(x_1,\cdots,x_d) dx_1=  \|u_1^{(1)}\|_{L^2(\Omega_1)}^2F^{(1)}(x_2,\ldots,x_d) + \int_{\Omega_1} u_1^{(1)}(x_1) G^{(1)}(x_1,\ldots,x_d)dx_1, 
\end{equation}
which implies that
$$ \int_{\Omega_1} u_1^{(1)}(x_1) G^{(1)}(x_1,\ldots,x_d)dx_1=0.$$

We proceed then to the next step of the inner iteration in the method, in which we take $F^{(1)}$ as the tensor to be approximated. We have suppose that the largest singular value corresponds again to the first of the remaining unfoldings and hence:

$$F^{(2)}(x_3,\ldots,x_d)=\int_{\Omega_{2}} F^{(1)}(x_2,\ldots,x_d)u_1^{(2,1)}(x_2)dx_2,$$

$$F^{(1)}=u_1^{(2,1)}(x_2)F^{(2)}(x_3,\ldots,x_d)+G^{(2)}(x_2,\ldots,x_d).$$

From which we can deduce:

$$\int_{\Omega_{2}} u_1^{(2,1)}(x_2) G^{(2)}(x_2,\ldots,x_d)dx_2=0.$$ \\

This results could be extended to all the substeps of the inner iteration, what means that for all $n \in \mathbb{N}$ from $1$ to $d$:\\
\begin{equation}F^{(n)}=F^{(n+1)}u_1^{(n+1,1)}+G^{(n+1)}.\end{equation}

And, for all $n \in \mathbb{N}$ from $1$ to $d$:\\
\begin{equation}\int_{\Omega_{n}} G^{(n)} u_1^{(n,n-1)} dx_n= 0 
\label{eq:orthGu}
\end{equation}

Let us consider the first iteration for which the residual $R^{(0)} = F-F^{(0)}=F$ and we can write:

$$R^{(1)}=F-c_1 u_1^{(1,1)}\otimes u_1^{(2,1)}\otimes...\otimes u_1^{(d-1,d-2)}\otimes u_1^{(d,d-1)},$$

where, in the following, to continue with the stablished notation: $ u_1^{(1,1)} = u_1^{(1)}$.

Using the above obtained relations this can be rewritten as:

\begin{equation*}
R^{(1)}= u_1^{(1,1)} u_1^{(2,1)} \ldots u_1^{(d-3,d-4)} u_1^{(d-2,d-3)} G^{(d-1)} + u_1^{(1,1)}... u_1^{(d-3,d-4)} G^{(d-2)} + \ldots + u_1^{(1,1)} G^{(2)} + G^{(1)}.
\end{equation*}

The squared norm of the residual reads:

$$
\|R^{(1)}\|_{L^2(\Omega)}^2= \int_\Omega ({u_1^{(1,1)}}^2 {u_1^{(2,1)}}^2\ldots {u_1^{(d-3,1)}}^2 {u_1^{(d-2,1)}}^2 {G^{(d-1)}}^2 +{u_1^{(1,1)}}^2 {u_1^{(2,1)}}^2 \ldots {u_1^{(d-3,1)}}^2 {u_1^{(d-2,1)}}^2 G^{(d-1)} G^{(d-2)} + \ldots$$
$$
\ldots + {u_1^{(1,1)}}^2 {u_1^{(2,1)}}^2 \ldots {u_1^{(d-3,1)}}^2 {G^{(d-2)}}^2 + \ldots +{u_1^{(1,1)}}^2 {G^{(2)}}^2 + {G^{(1)}}^2)dx_1 \ldots dx_d
$$\\

Separating the integrals, applying the orthogonality relation obtained in Eq.\eqref{eq:orthGu} and using the fact that the modes are unitary, the crossed terms in $G$ vanish and this simplifies as:

$$
\|R^{(1)}\|_{L^2(\Omega)}^2= \int_\Omega ({u_1^{(1,1)}}^2 {u_1^{(2,1)}}^2 \ldots {u_1^{(d-3,d-4)}}^2 {u_1^{(d-2,d-3)}}^2 {G^{(d-1)}}^2 +{u_1^{(1,1)}}^2 {u_1^{(2,1)}}^2 \ldots {u_1^{(d-3,d-4)}}^2 {G^{(d-2)}}^2 + \ldots $$

$$ \ldots +{u_1^{(1,1)}}^2 {G^{(2)}}^2 +{G^{(1)}}^2)dx_1 \ldots dx_d
$$

Replacing the norms by the squared singular values due to the properties of $G$:
\begin{equation}
 \int_{\Omega_k \times \ldots \times \Omega_d} G^{(k)} dx_k...dx_{d-1}dx_d= \sum_{j>1} {\sigma_j^{(k, )}}^2 
 \label{eq:Gprop}
\end{equation}

using this result in the expression of the residual is possible to get the analogue expression for the total residual as the sum of the reminders on each term.

$$\|R^{(1)}\|_{L^2(\Omega)}^2=\|G^{(1)}\|_{L^2(\Omega)}^2+\|G^{(2)}\|_{L^2(\Omega_2 \times \ldots \times \Omega_d)}^2+...+ \|G^{(d-1)}\|_{L^2(\Omega_{d-1} \times \Omega_d)}^2$$\\

\paragraph{Example:Properties in the $4$-th dimensional case}
In the following we illustrate these steps on a $4-$th order tensor. After the first iteration, the residual reads:
\begin{equation*}
R^{(1)}=F-c_1 u_1^{(1,1)}\otimes u_1^{(2,1)}\otimes u_1^{(3,2)} \otimes u_1^{(4,3)}= u_1^{(1,1)}u_1^{(2,1)}G^{(3,2)} + u_1^{(1,1)}G^{(2)}+ G^{(1)}
\end{equation*}

Where we have taken into account the fact that in the last iteration $F^{(3)}=c_1u_1^{(4,1)}$ because of its orthogonality in the POD decomposition.\\

The squared norm of the residual reads:
$$
\|R^{(1)}\|_{L^2(\Omega)}^2= \int_{\Omega}({u_1^{(1,1)}}^2 {u_1^{(2,1)}}^2 {G^{(3)}}^2 + 2 {u_1^{(1,1)}}^2 u_1^{(2,1)} G^{(3)}G^{(2)} + 2 {u_1^{(1,1)}}^2 G^{(2)} G^{(1)} + {u_1^{(1,1)}}^2 {G^{(2)}}^2 + \ldots $$
$$ \ldots + 2u_1^{(1,1)} G^{(2)} G^{(1)}  + {G^{(1)}}^2 )dx_1\ldots dx_4 $$
Using the orthogonality relations seen in Eq.\eqref{eq:orthGu} and the orthonormality of the modes:
\begin{equation*}
\|R^{(1)}\|_{L^2(\Omega)}^2= \int_{\Omega_{3} \times \Omega_4}  {G^{(3)}}^2 dx_3 dx_4 + \int _{\Omega_2 \times \Omega_{3} \times \Omega_4}{G^{(2)}}^2 dx_2 dx_3 dx_4 + \int_{\Omega}  {G^{(1)}}^2 dx_1 dx_2 dx_3 dx_4
\end{equation*}\\

Adding the proprieties of $G$ from Eq.\eqref{eq:Gprop} it leads to:

\begin{equation*}
\|R^{(1)}\|_{L^2(\Omega)}^2= {\sigma_1^{(3,2)}}^2 -{\sigma_1^{(4,3)}}^2 +{\sigma_1^{(2,1)}}^2-{\sigma_1^{(2,1)}}^2+\|F\|_{L^2(\Omega)}^2-{\sigma_1^{(3,2)}}^2 =\|F\|_{L^2(\Omega)}^2-c^2
\end{equation*}

where we have called $c$ to the last of the considered eigenvalues.\\

What we can see at this point is that $\|F^{(1)}\|_{L^2(\Omega)}^2={\sigma_1^{(1)}}^2$.\\

And
$$\|F\|_{L^2(\Omega)}^2=\|F^{(1)}\|_{L^2(\Omega)}^2+\|G^{(1)}\|_{L^2(\Omega)}^2 $$

We can also confirm that the squared norm of the total residual in one iteration is the sum of the squared norms of the residuals in each term on the iteration:

$$\|R^{(1)}\|_{L^2(\Omega)}^2=\|G^{(1)}\|_{L^2(\Omega)}^2+\|G^{(2)}\|_{L^2(\Omega_{2} \times \Omega_3 \times \Omega_4)}^2+ \|G^{(3)}\|_{L^2(\Omega_{3} \times \Omega_4)}^2$$

\subsection{CP-TT is optimal for $d=2$, and it retrieves exactly a rank$-1$ tensor}
The method proposed, when $d=2$, reduces to compute a POD of a two variables function $F(x_1,x_2)$. By the Eckart-Young theorem, the result is optimal. As for other methods, when the number of variables is $d>2$, CP-TT is in general sub-optimal. 
\medskip

When the function to be approximated is a rank$-1$ tensor, it is exactly recovered by the CP-TT method in $N_{it}=1$ iteration. 

Let $F = f^{(1)}\otimes\ldots \otimes f^{(d)}$. Without loss of generality, let the order of the chosen best unfolding be $I=\left\lbrace 1,\ldots, d\right\rbrace$. 
The best approximation for the first unfolding gives us the term $u^{(1)}_1 = \frac{f^{(1)}}{\| f^{(1)} \|_{L^2(\Omega^{(1)})}}$. Proceeding with the CP-TT iteration, and considering the unfolding of $F^{(1)}$ relative to $x_2$, the result is $u^{(2,1)} = \frac{f^{(2)}}{\| f^{(2)} \|_{L^2(\Omega^{(2)})}}$. This can be iterated, leading to the following outcome of the first CP-TT iteration:
\begin{eqnarray}
\tilde{F} = c u^{(1)}\otimes \ldots \otimes u^{(d)}, \nonumber \\
c = \| f^{(1)} \|_{L^2(\Omega^{(1)})}\cdot \ldots \cdot \| f^{(d)} \|_{L^2(\Omega^{(d)})}, \nonumber \\
u^{(1)} = \frac{f^{(1)}}{\| f^{(1)} \|_{L^2(\Omega^{(1)})}}, \nonumber \\
\ldots \nonumber \\
u^{(d)} = \frac{f^{(d)}}{\| f^{(d)} \|_{L^2(\Omega^{(d)})}} \nonumber .
\end{eqnarray}
Henceforth, $\tilde{F}=F$, and the method stops.

\subsection{Computing $k$ terms: stability}
The CP-TT iteration can be suitably modified to compute more than a rank-one update. Let $k>1$ be the number of pure tensor terms to be computed. We introduce the following modification of the above described method. When computing the unfoldings POD for $M_i$, instead of selecting the index of the variable corresponding to the largest first singular value, solve the following problem:
\begin{equation}
\label{eq:modified_opt}
i_* = \mathrm{arg}\max_{1\leq i\leq d} \sum_{j=1}^k \left(\sigma_j^{(i)}\right)^2.
\end{equation}
Then, consider the set of $u^{(i_*)}_j$, $j=1,\ldots,k$. These can be used to compute:
\begin{eqnarray}
F^{(1)}_j = \int_{\Omega_{i_*}} F u^{(i_*)}_j \ dx_{i_*}. \nonumber
\end{eqnarray}
The CP-TT iteration can be carried out independently for each of the $F^{(1)}_j$, for the remaining variables, leading to the definition of $k$ pure tensor terms. The method so defined computes a stable CP decomposition, in the sense that:
\begin{equation}
\| \sum_{j=1}^k t_j \|_{L^2(\Omega)}^2 < C \ \Rightarrow \ \| t_j \|_{L^2(\Omega)}^2 < C_j, \ \forall j. \nonumber
\end{equation}
This follows from the very first step of the modified iteration. Indeed:
\begin{equation}
\| \sum_{j=1}^k t_j \|_{L^2(\Omega)}^2 = \sum_{l=1}^k\sum_{m=1}^k \langle t_l, t_m\rangle = \sum_{l,m=1}^k \langle u^{(i_*)}_l,u^{(i_*)}_m \rangle_{L^2(\Omega_{i_*})} \langle w_l^{(i_*)},w_m^{(i_*)}\rangle_{ L^2(\Omega\setminus\Omega_{i_*})  }
\end{equation}
Since the terms $u^{(i_*)}_j$ are elements of an orthonormal basis of $L^2(\Omega_{i_*})$ we have:
\begin{equation}
\sum_{l,m=1}^k \langle u^{(i_*)}_l,u^{(i_*)}_m \rangle_{L^2(\Omega_{i_*})} \langle w_l^{(i_*)},w_m^{(i_*)}\rangle_{ L^2(\Omega\setminus\Omega_{i_*})  } = \sum_{l=1}^k \| t_l \|^2_{L^2(\Omega)},
\end{equation}
and this is enough to conclude.

\subsection{Discretization and computational costs}
In this section, the computational cost of the proposed method is detailed. Let $(\mathcal{N}_1,\ldots, \mathcal{N}_d)$ be the number of degrees of freedom used to discretize the problem in every direction. The storage of a full tensor would require a memory of $\prod_{j=1}^d \mathcal{N}_j$. The simplicity of the CP tensor format is also complemented by its appealing storage scaling, which is: $r \cdot \sum_{j=1}^d \mathcal{N}_j$.
In view of implementing the proposed method, some aspects need to be considered. The most costly operation to be performed is the computation of the POD of the $d$ unfoldings. In the following, a method is proposed, in the case in which the function to be compressed is given in CP format, to provide a cheap alternative for the POD decomposition of the unfoldings.

In \cite{badeau2008fast} a fast multilinear singular value decomposition is proposed for symmetric Toeplitz and Hankel third order tensors has been proposed. A fast computation of the gradients in ALS method is proposed in \cite{phan2013fast}. In \cite{vervliet2019exploiting} efficient representations are exploited to speed up large-scale tensor decompositions. In \cite{kaya2015scalable} the authors propose parallel methods to accelerate the decomposition of sparse tensors.

\subsubsection{Computing the POD of an unfolding without storing it}
Let the tensor $F$ be given in CP format:
$$F = \sum_{i=1}^r c_i f_i^{(1)}\otimes \ldots \otimes f_i^{(d)},$$
where $\forall i, \ f_i^{(j)}\in L^2(\Omega_j)$ and $c_i\in\mathbb{R}$. Henceforth, $rank(F)=r$.

Without loss of generality we present the method for $M^{(1)}(x_1,y_1)$, the unfolding relative to the first variable.

The first step consists in introducing the correlation operator:
\begin{eqnarray}
\forall f(x_1) \in L^2(\Omega_1), \ \ K^{(1)} f = \int_{\Omega_1} \kappa^{(1)}(x_1,x_1') f(x_1') \ dx_1',
\end{eqnarray}
where the kernel $\kappa^{(1)}$ is defined as follows:
\begin{equation}
\kappa^{(1)}(x_1,x_1') = \int_{\Omega_2\times\ldots\times\Omega_d} M(x_1,y_1)M(x_1',y_1) \ dy_1.
\end{equation}
The operator $K^{(1)}$ is compact, non-negative, self-adjoint, and hence, there exists a sequence of orthonormal functions $u^{(1)}_i(x_1), \ i\in\mathbb{N}^*$ and a sequence of scalar values $\lambda_i\geq 0$ such that:
\begin{equation}
K^{(1)} u^{(1)}_i = \lambda_i u^{(1)}_i.
\end{equation}

When the tensor is given in CP format, the expression of the kernel is reduced to:
\begin{equation}
\kappa^{(1)}(x_1,x_1') = \sum_{i,j=1}^r c_i c_j f_i^{(1)}(x_1)f_j^{1}(x_1') \int_{\Omega_2} f_i^{(2)}(x_2)f_j^{2}(x_2) \ dx_2 \ldots \int_{\Omega_d} f_i^{(d)}(x_d)f_j^{d}(x_d)dx_d.
\end{equation}

The POD decomposition of the fibers $f_i^{(1)}$ is computed, providing:
\begin{equation}
\label{eq:PODfibers}
f_i^{(1)}(x_1) = \sum_{m=1}^r a_{im} z^{(1)}_m(x_1),
\end{equation}
where $\langle z^{(1)}_m, z^{(1)}_l \rangle_{L^2(\Omega_1)}=\delta_{lm}$ and $a_{im} = \sum_{s=1}^r\sigma^{(1)}_s \delta_{is} v_{sm}$, where $\sigma^{(1)}_s\geq 0$ and $v_{sm}$ are the entries of a unitary matrix. This is injected into the correlation operator, leading to:
\begin{eqnarray}
\kappa^{(1)}(x_1,x_1') = \sum_{l,m=1}^r z^{(1)}_l(x_1) z^{(1)}_m(x_1')\left[\sum_{i,j=1}^r a_{il}a_{jm} c_i c_j \prod_{k=2}^d \int_{\Omega_k} f_i^{(k)}(x_k)f_j^{k}(x_k) \ dx_k \right].
\end{eqnarray}
Let $A\in\mathbb{R}^{r\times r}$ be the matrix defined as:
\begin{equation}
\label{eq:matrixA}
A_{lm} = \left[\sum_{i,j=1}^r a_{il}a_{jm} c_i c_j \prod_{k=2}^d \int_{\Omega_k} f_i^{(k)}(x_k)f_j^{k}(x_k) \ dx_k \right].
\end{equation}
This matrix is symmetric and positive semidefinite by construction. Henceforth, there exists a complete orthonormal basis of $\mathbb{R}^r$ that diagonalises it: $A = W S W^T$. The eigenfunctions of the correlation operator, which are also the left POD modes of the first unfolding are defined as:
\begin{equation}
u^{(1)}_i(x_1) = \sum_{l=1}^r z_l(x_1) W_{li}.
\end{equation}
The associated singular value is $\sigma_i^{(1)} = \sqrt{S_{ii}}$.

The method is summarised in Algorithm \ref{alg:UnU}.

\begin{algorithm}[]
\caption{POD decomposition of the unfolding}
\begin{algorithmic}
 \label{alg:UnU}
\STATE{\textbf{Input:}}
\STATE{ $M^{(1)}(x_1,y_1)$}
\medskip
\STATE{\textbf{Method:}}
\STATE{Compute POD of the fibers $f_i^{(1)}$, see Eq.\eqref{eq:PODfibers}}
\STATE{Assemble the matrix A, see Eq.\eqref{eq:matrixA}}
\STATE{Diagonalise it: $A = W S W^T$}
\STATE{Compute the modes $u_i^{(1)}$ and the associated $\sigma_i^{(1)}$}
 \end{algorithmic}
\end{algorithm}

\subsubsection{Computational cost} 
The computational cost of the method is detailed hereafter. 

We start by describing the cost of the computation of the POD for one given unfolding, which is the most expensive operation of the CP-TT iteration. There are two distinct cases. If $\mathcal{N}_i<r$, then, the discrete counterpart of $\kappa^{(i)}$ can be evaluated, and its eigenvalue decomposition computed directly. This leads to a cost of the form:
\begin{itemize}
\item Correlation assembly: $\mathcal{O}\left( r^2 (\sum_{j\neq i}^d c_I \mathcal{N}_j) + (2r-1)(r\mathcal{N}_i + \mathcal{N}_i^2) \right)$.
\item Eigenvalue decomposition: $\mathcal{O}\left( \mathcal{N}_i^3\right)$.
\end{itemize}

If $\mathcal{N}_i>r$, then, the method presented in the previous section is used and its cost reads:
\begin{itemize}
\item SVD of the fibers $f^{(i)}$: $\mathcal{O}\left(\mathcal{N}_i r^2\right)$.
\item Assembly of $A$: $\mathcal{O}\left( r^2 (\sum_{j\neq i}^d c_I \mathcal{N}_j)  + r^2(1+2(2r-1)) + r^3\right)$.
\item Eigenvalue decomposition of $A$: $\mathcal{O}\left(r^3\right)$.
\item Computation of the left modes: $\mathcal{O}\left( (2r-1)\mathcal{N}_i \right)$.
\end{itemize}

Overall, the cost of this stage is of order $\mathcal{O}\left(r^2 d \mathcal{N} + \min \left\lbrace \mathcal{N}_i^3, \ r^3 \right\rbrace \right)$.
 
\medskip

When considering the CP-TT iteration, at the first step we compute $d$ POD decompositions (these could be performed in a parallel way). After this stage, the tensor is tested against the mode corresponding to the chosen direction. The cost is $\mathcal{O}\left( r \mathcal{N}_i \right)$. Then, the POD of $d-1$ unfoldings is computed. Remark that the cost is the same as the one detailed above, with the exception that the operations involving the assembly of the matrices costs less (since we are working at dimension $d-1$). 

\paragraph*{Remark}
The first stage of the computation is similar, in a way, to the HOSVD method (\cite{de2000multilinear}). However, in the present method, only the first singular triplet is needed and actually used, which leads to a less expensive computation. Moreover, no storage of the core tensor is required.

\section{Numerical Experiments}
\label{sec:numExp}

Once we have studied the properties of the method theoretically, we will proceed to do some numerical experiments based on function compression. In this section we can get an idea of how the method acts in practice. The proposed tests are made for functions which admit a Fourier decomposition involving a finite number of modes. In the tests we do a comparison between three methods, namely: ALS, ASVD and CPTT. \\

Let $(x_1, \ldots, x_d) \in \Omega = [0,1]^d$.
Let $(k_1,...,k_d)\in\mathbb{N}^d$ be the wave numbers. The function to be compressed is assumed to be given in CP format :
\begin{equation}
\label{refF}
F(x_1,\ldots,x_d) = \sum_{k_1=1}^{l_1}\sum_{k_2=1}^{l_2}\ldots\sum_{i_d=1}^{l_d} a_{k_1\ldots k_d} \sin(\pi k_1 x_1)\times\ldots\times\sin(\pi k_d x_d)
\end{equation}

Let $\beta>0$. The values of $\left(  l_i\right)_{1\leq i\leq d}$ are chosen to be a family of independent random integers uniformly distributed between 1 and 6. Let  $\left(\alpha_{k_1\ldots k_d}\right)_{k_1,...,k_d}$ be a family of independent random variables uniformly distributed in $ [-1,1]$. The amplitudes $a_{k_1\ldots k_d}$ are defined as:
$$ a_{k_1\ldots k_d}=\frac{\alpha_{k_1\ldots k_d}}{{(\sqrt{k_1^2 + ...+ k_d^2})}^\beta}.$$

For different random samples we obtain different functions $F$ with the form presented in Eq.\eqref{refF} , the amplitude change and we obtain different functions preserving the shape of $F$. The value of the parameter $\beta$ determines the regularity of the functions in the sense that for $s\in\mathbb{N}$ it holds:
\begin{equation}
\beta > \frac{d}{2} + s \ \Rightarrow \ F \in H^{s}(\Omega).
\end{equation}

We are testing how the three methods behave for the compression of 32 different functions generated by the random procedure described above for values of $d$ ranging from 4 to 16. ALS and ASVD are both fixed point based methods, the tolerance for the fixed point has been set as $1.0 \times 10^{-4}$ and the maximum number of iterations of the method $it_{max}=100$.  A uniform discretization grid of $\Omega$ with 25 degrees of freedom per direction is used for the discretization of $F$.\\

\subsection{Results for functions in $L^2(\Omega)$} 

\begin{figure}[h!]
\includegraphics[width=7cm]{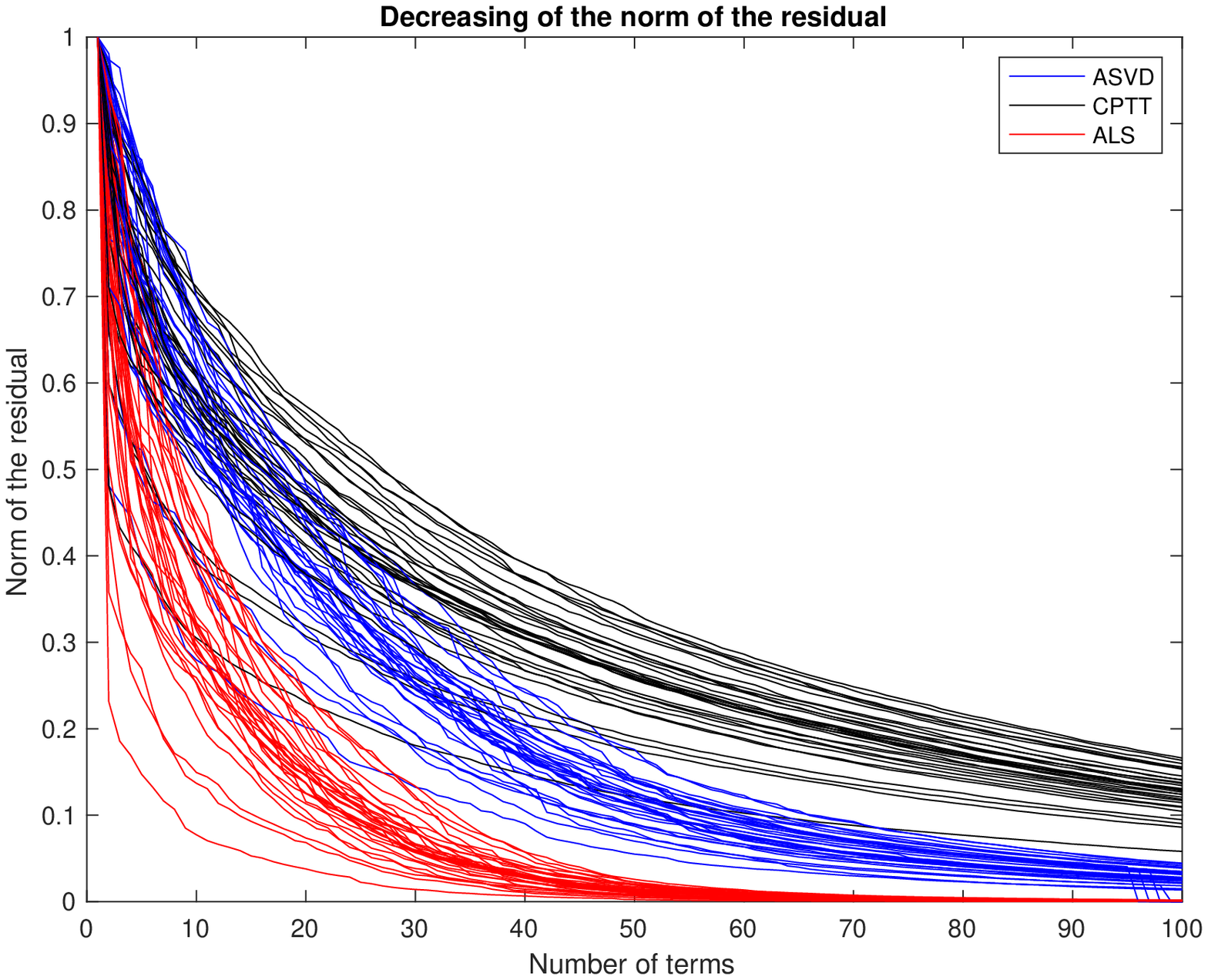}
\includegraphics[width=7cm]{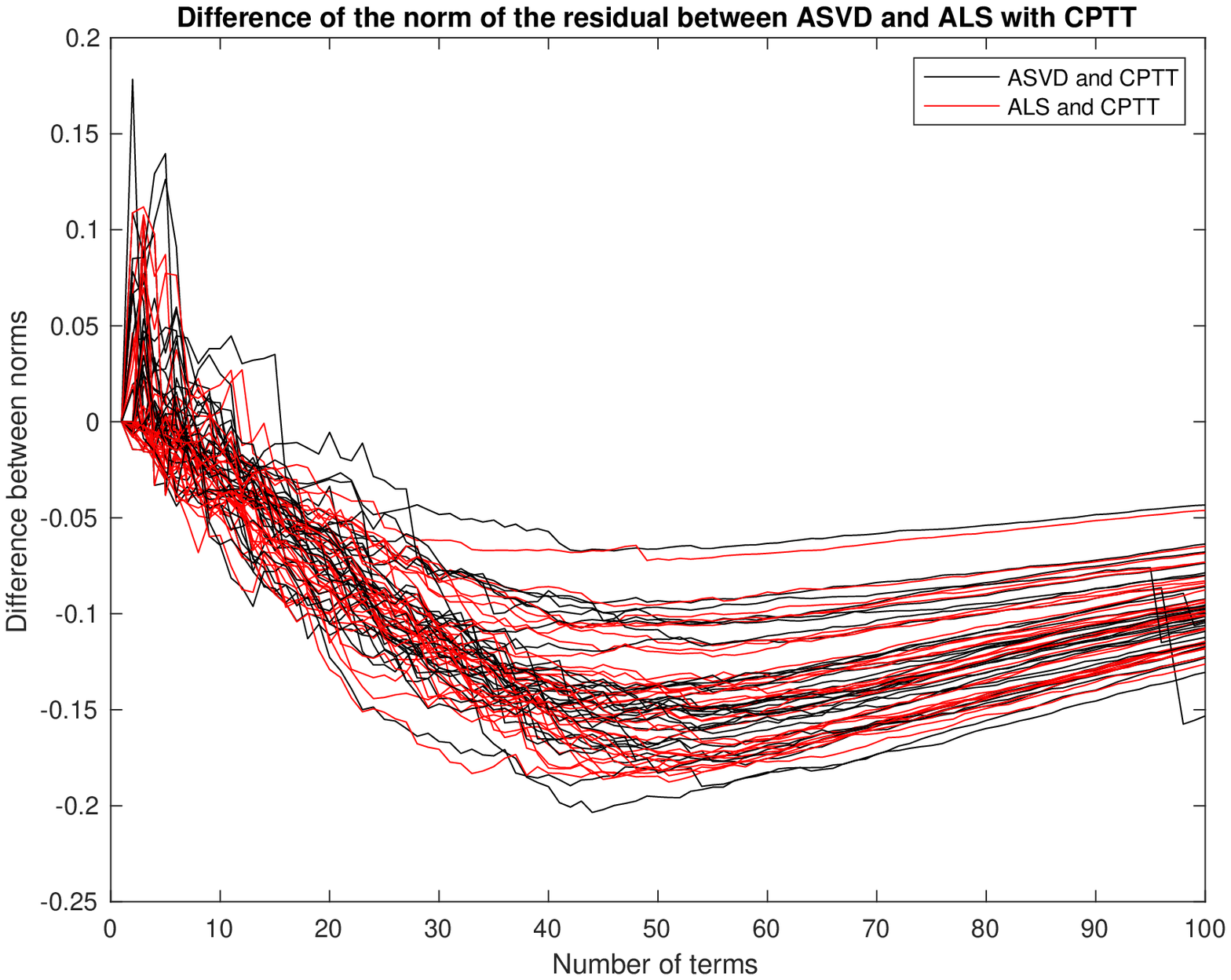}
\caption{Plot associated to functions described in Eq.\eqref{refF}. On the left, the decreasing of the $L^2$ norm of the difference between the exact function $F$ and its approximation in the 4-dimensional tensor of the 32 $L^2$ functions, $\beta=2.1$, using ALS in red, ASVD in blue and the CP-TT method in black. On the right, the difference between both norms the difference between the exact function $F$ and its approximation $\|R\|_{ALS}-\|R\|_{CP-TT}$ in red and $\|R\|_{ASVD}-\|R\|_{CP-TT}$ in black.}
\label{refI2}
\end{figure}

\begin{figure}[h!]
\includegraphics[width=7cm]{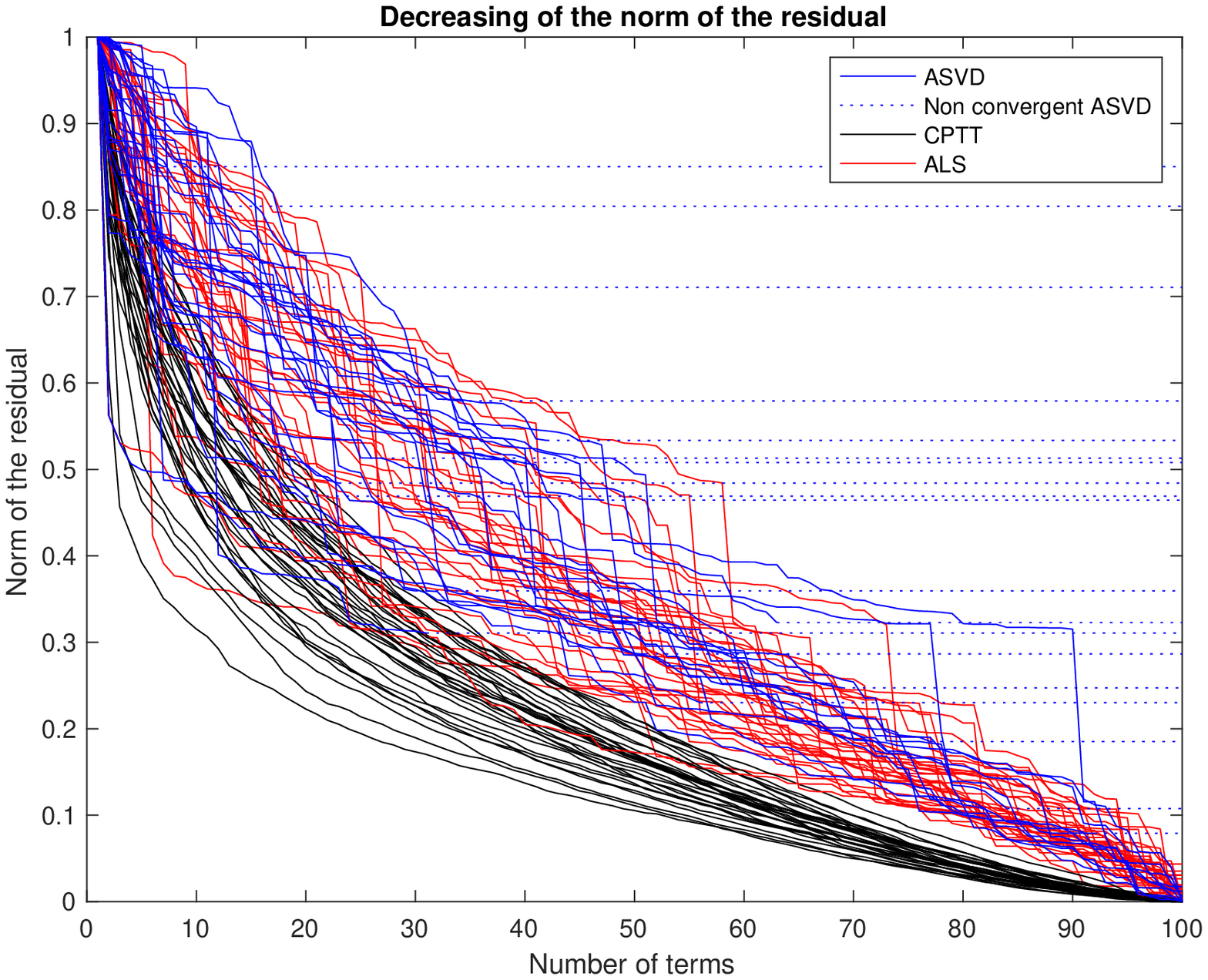}
\includegraphics[width=7cm]{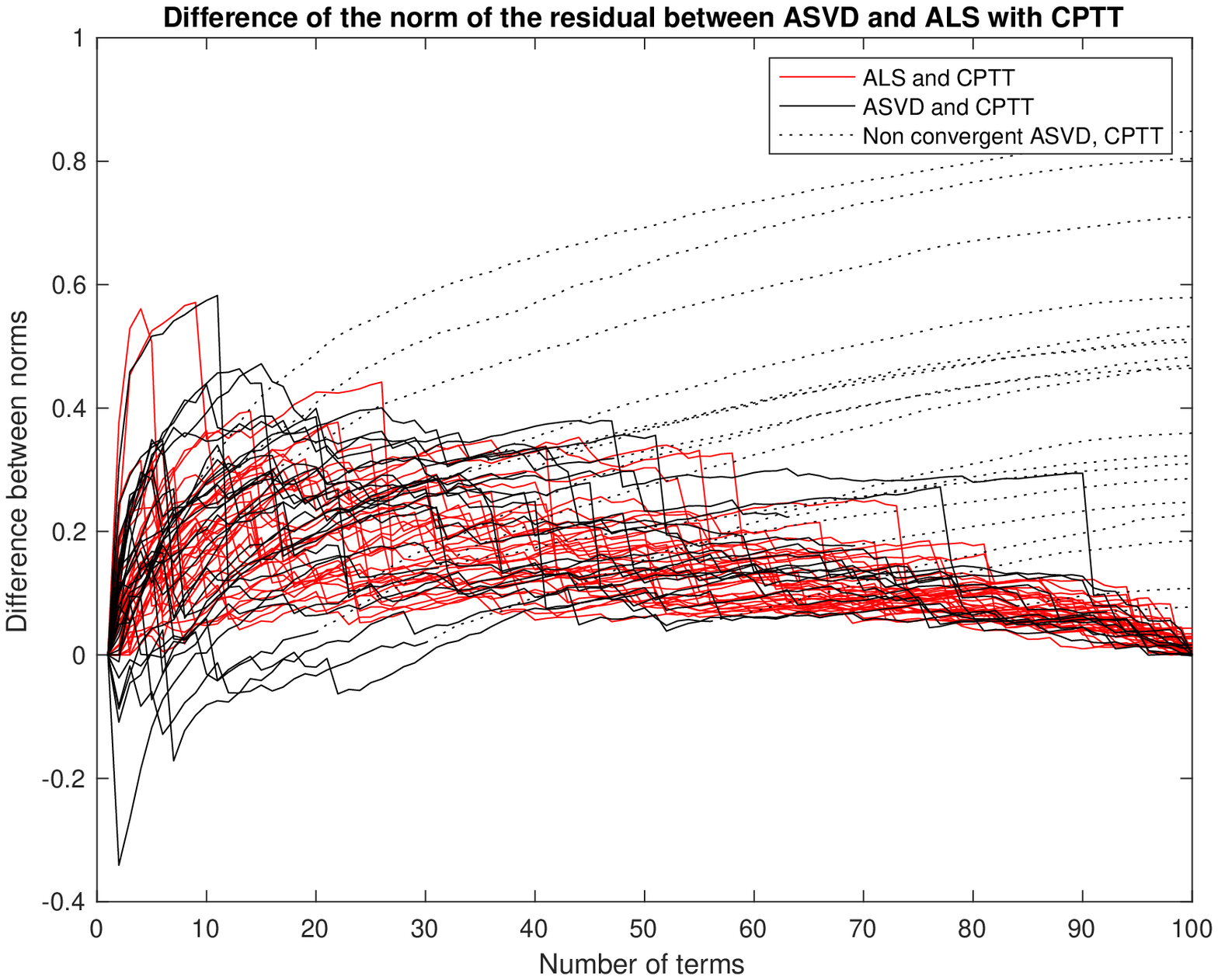}
\caption{Plot associated to functions described in Eq.\eqref{refF}. On the left, the decreasing of the $L^2$ norm of the difference between the exact function $F$ and its approximation in the 12-dimensional tensor of the 32 $L^2$ functions, $\beta=6.1$, using ALS in red, ASVD in blue and the CP-TT method in black. The dotted line shows as a constant the last point that converged on the ASVD method. On the right, the difference between both norms the difference between the exact function $F$ and its approximation $\|R\|_{ALS}-\|R\|_{CP-TT}$ in red and $\|R\|_{ASVD}-\|R\|_{CP-TT}$ in black.}
\label{refI3}
\end{figure}

\begin{figure}[h!]
\includegraphics[width=7cm]{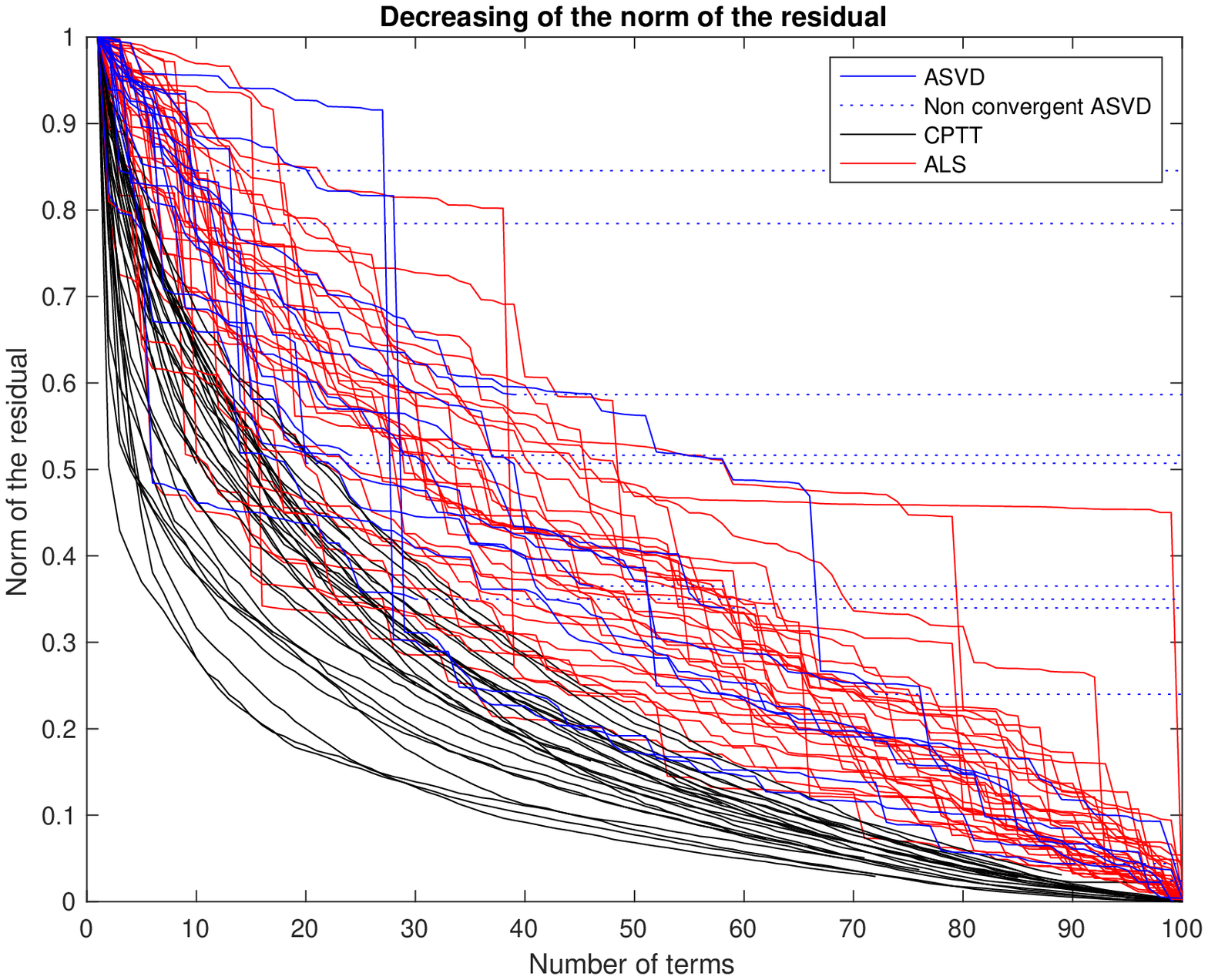}
\includegraphics[width=7cm]{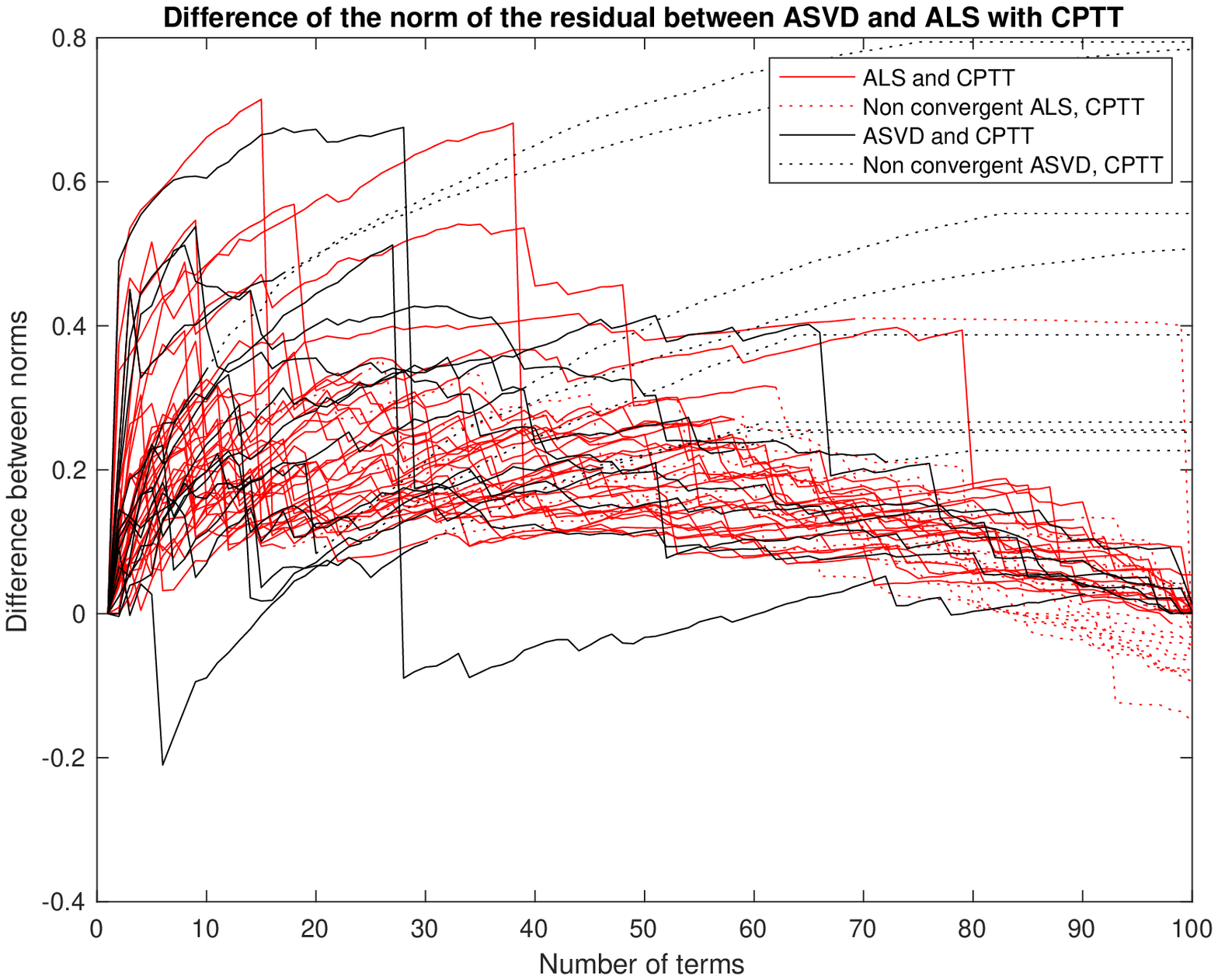}
\caption{Plot associated to functions described in Eq.\eqref{refF}. On the left, the decreasing of the $L^2$ norm of the difference between the exact function $F$ and its approximation in the 16-dimensional tensor of the 32 $L^2$ functions, $\beta=8.1$, using ALS in red, ASVD in blue and the CP-TT method in black. The dotted line shows as a constant the last point that converged on the ASVD and CP-TT methods. On the right, the difference between both norms the difference between the exact function $F$ and its approximation $\|R\|_{ALS}-\|R\|_{CP-TT}$ in red and $\|R\|_{ASVD}-\|R\|_{CP-TT}$ in black.}
\label{refI5}
\end{figure}

Firstly, the method is applied to some test cases in which the functions belong to $L^2(\Omega)$, namely when the value of the parameter $\beta$ is chosen to be equal to $\frac{d}{2}+0.1$. On the left hand side of Figure \ref{refI2} (respectively Figure \ref{refI3} and Figure \ref{refI5}) , the $L^2$ norm of the difference between the exact function $F$ and its approximation computed by one of the three algorithms is plotted as a function of the rank of the approximation, where $d=4$ (respectively $d=12$ and $d=16$).

Whereas at dimension 4 (Figure \ref{refI2}) ALS outperforms both ASVD and CP-TT, at dimensions $d=12$ and $d=16$ (figures \ref{refI3} and \ref{refI5} respectively) CP-TT is featured by a better behavior. In particular, the compression rate is better on average and the decrease of the norm of the error with respect to the rank of the approximation is more regular. To better highlight this, on the right-hand side of figures \ref{refI2}, \ref{refI3} and \ref{refI5} we plot the difference of the error norms of ALS and CP-TT and ASVD and CP-TT as function of the rank of the approximation.\\

\begin{figure}[h!]
\includegraphics[width=7cm]{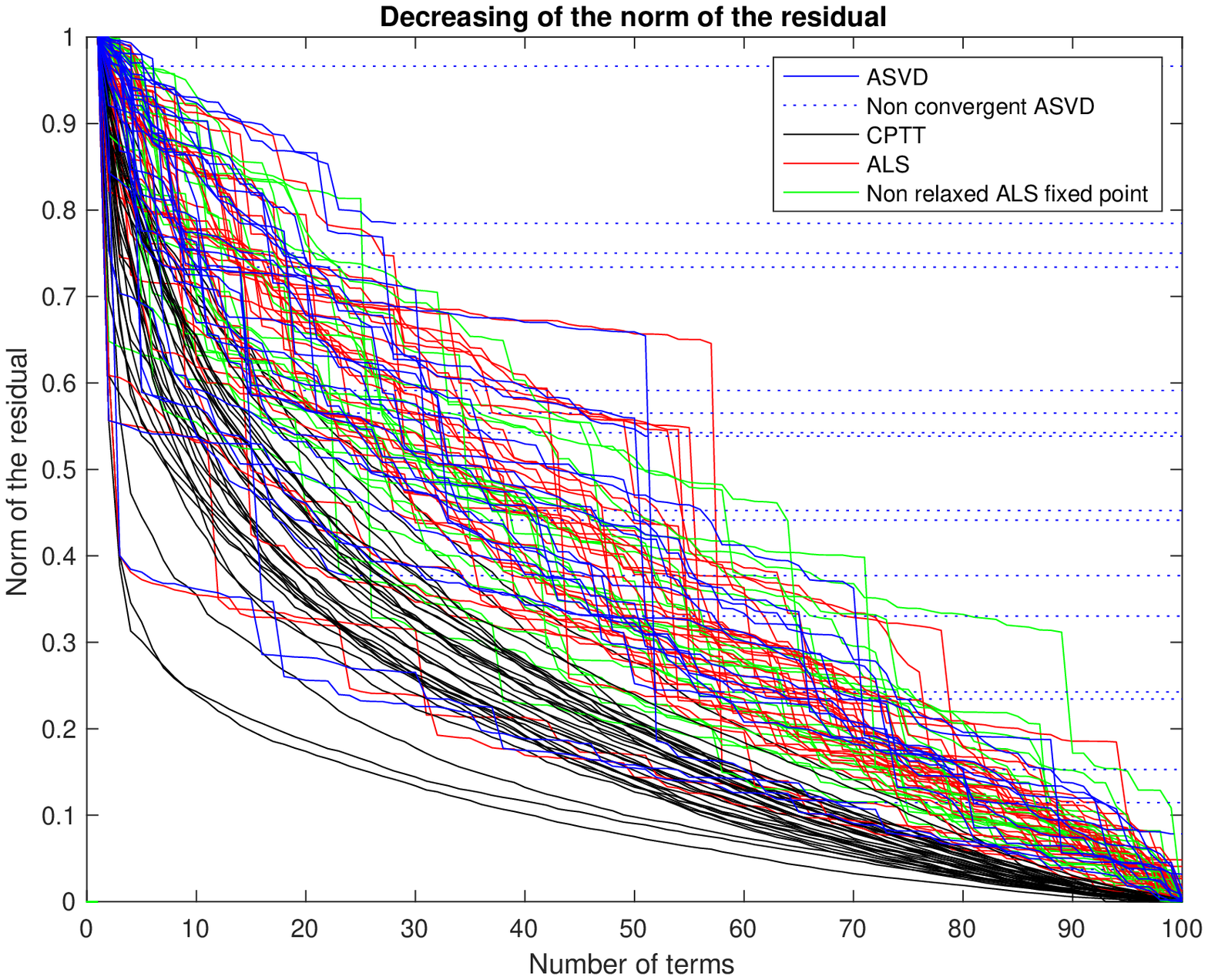}
\includegraphics[width=7cm]{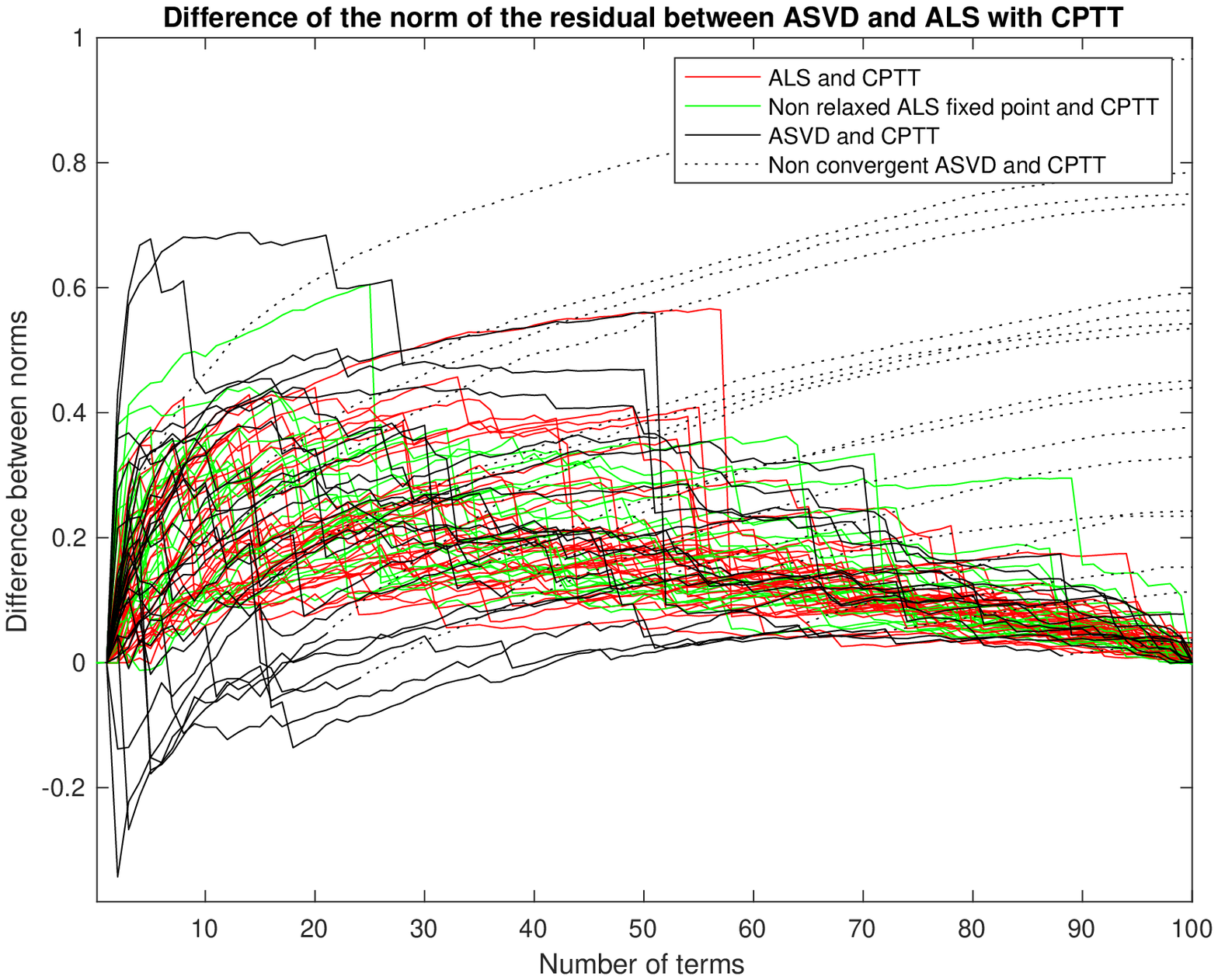}
\caption{Plot associated to functions described in Eq.\eqref{refF}. On the left, the decreasing of the $L^2$ norm of the difference between the exact function $F$ and its approximation in the 12-dimensional tensor 
of the 32 $L^2$ functions, $\beta=6.1$, using ALS in red, ALS without the relaxed fixed point in green, the ASVD in blue and the CP-TT method in black. The dotted line shows as a constant the last point that converged on the ALS method. On the right, the difference between both norms the difference between the exact function $F$ and its approximation $\|R\|_{ALS}-\|R\|_{CP-TT}$ in red and $\|R\|_{ASVD}-\|R\|_{CP-TT}$ in black.}
\label{refI4}
\end{figure}

In the previous test a relaxed version of ALS (standard ALS with more flexibility in the scalar products in order to optimize the coefficients) was used that is observed to converge more often than the standard version of ALS. In Figure \ref{refI4},  the relaxed and the non relaxed version of the ALS algorithm in the 12-dimensional case are compared. Both versions were observed to yield similar compression behavior.\\

The results for $L^2(\Omega)$ functions, when changing the dimension, are reported in \ref{tab:tableL2}. In particular, we show the mean and the standard deviation of the error (on the 32 random functions) when the approximation rank is $(25,50,75)$ for ALS, ASVD and CP-TT. From the results, we can see that for low-order tensors (\emph{e.g.} when $d=4$) ALS has better performances, whereas for higher order tensors CP-TT outperforms the other methods both in terms of mean and standard deviation (suggesting that it somehow enjoys a certain stability).

\begin{table}[tbhp]
\begin{center}
\begin{tabular}{cc | ccc | ccc | }
\cline{3-8}
                                &                                      & \multicolumn{1}{c}{}     & \multicolumn{1}{c}{Mean}  & \multicolumn{1}{c|}{}     & \multicolumn{1}{c}{}     & \multicolumn{1}{c}{Standard deviation} & \multicolumn{1}{c|}{}     \\ \hline
\multicolumn{1}{|c|}{Dimension ($d$)} & \multicolumn{1}{c|}{Rank ($r$)} & \multicolumn{1}{c|}{ALS} & \multicolumn{1}{c|}{CPTT} & \multicolumn{1}{c|}{ASVD} & \multicolumn{1}{c|}{ALS} & \multicolumn{1}{c|}{CPTT}              & \multicolumn{1}{c|}{ASVD} \\ \hline
\multicolumn{1}{|l|}{}          & 25                                   & \multicolumn{1}{c|}{0.2942}   & \multicolumn{1}{c|}{0.3826}    & 0.3118                         & \multicolumn{1}{c|}{0.0702}    & \multicolumn{1}{c|}{0.0850}                  &0.0843                           \\
\multicolumn{1}{|c|}{4}         & 50                                   & \multicolumn{1}{c|}{0.1082}    & \multicolumn{1}{c|}{0.2433}     &0.1257                           & \multicolumn{1}{c|}{0.0326}    & \multicolumn{1}{c|}{0.0568}                  &0.0664                           \\
\multicolumn{1}{|l|}{}          & 75                                   & \multicolumn{1}{c|}{0.0508}    & \multicolumn{1}{c|}{0.1681}     &0.0689                           & \multicolumn{1}{c|}{0.0180}    & \multicolumn{1}{c|}{0.0408}                  &0.0666                           \\ \hline

\multicolumn{1}{|l|}{}          & 25                                   & \multicolumn{1}{c|}{0.4479}   & \multicolumn{1}{c|}{0.3771}    & 0.4806                         & \multicolumn{1}{c|}{0.1099}    & \multicolumn{1}{c|}{0.0826}                  &0.1074                           \\
\multicolumn{1}{|c|}{6}         & 50                                   & \multicolumn{1}{c|}{0.2705}    & \multicolumn{1}{c|}{0.1982}     &0.2883                           & \multicolumn{1}{c|}{0.0752}    & \multicolumn{1}{c|}{0.0485}                  &0.0675                           \\
\multicolumn{1}{|l|}{}          & 75                                  & \multicolumn{1}{c|}{0.1232}    & \multicolumn{1}{c|}{0.0806}     &0.1369                           & \multicolumn{1}{c|}{0.0325}    & \multicolumn{1}{c|}{0.0252}                  &0.0368                           \\ \hline

\multicolumn{1}{|l|}{}          & 25                                   & \multicolumn{1}{c|}{0.5341}   & \multicolumn{1}{c|}{0.3707}    & 0.5532                         & \multicolumn{1}{c|}{0.1183}    & \multicolumn{1}{c|}{0.0592}                  &0.1238                           \\
\multicolumn{1}{|c|}{8}         & 50                                   & \multicolumn{1}{c|}{0.3060}    & \multicolumn{1}{c|}{0.1909}     &0.3415                           & \multicolumn{1}{c|}{0.0722}    & \multicolumn{1}{c|}{0.0341}                  &0.0932                           \\
\multicolumn{1}{|l|}{}          & 75                                   & \multicolumn{1}{c|}{0.1592}    & \multicolumn{1}{c|}{0.0682}     &0.1807                           & \multicolumn{1}{c|}{0.0435}    & \multicolumn{1}{c|}{0.0160}                  &0.0625                           \\ \hline

\multicolumn{1}{|l|}{}          & 25                                  & \multicolumn{1}{c|}{0.5023}   & \multicolumn{1}{c|}{0.3598}    &0.5451                       & \multicolumn{1}{c|}{0.0879}    & \multicolumn{1}{c|}{0.0643}                  &0.1055                           \\
\multicolumn{1}{|c|}{10}         & 50                                & \multicolumn{1}{c|}{0.3191}    & \multicolumn{1}{c|}{0.1826}     &0.3797                           & \multicolumn{1}{c|}{0.0643}    & \multicolumn{1}{c|}{0.0342}                  &0.0774                           \\
\multicolumn{1}{|l|}{}          & 75                                   & \multicolumn{1}{c|}{0.1714}    & \multicolumn{1}{c|}{0.0655}     &0.2792                           & \multicolumn{1}{c|}{0.0453}    & \multicolumn{1}{c|}{0.0162}                  &0.1265                           \\ \hline

\multicolumn{1}{|l|}{}          & 25                                   & \multicolumn{1}{c|}{0.5170}   & \multicolumn{1}{c|}{0.3246}    & 0.5639                          & \multicolumn{1}{c|}{0.1117}    & \multicolumn{1}{c|}{0.0576}                  &0.1250                           \\
\multicolumn{1}{|c|}{12}         & 50                               & \multicolumn{1}{c|}{0.3249}    & \multicolumn{1}{c|}{0.1623}     & 0.4206                          & \multicolumn{1}{c|}{0.0824}    & \multicolumn{1}{c|}{0.0286}                  &0.1579                           \\
\multicolumn{1}{|l|}{}          & 75                                   & \multicolumn{1}{c|}{0.1543}    & \multicolumn{1}{c|}{0.0579}     &0.3498                           & \multicolumn{1}{c|}{0.0369}    & \multicolumn{1}{c|}{0.0113}                  &0.2057                           \\ \hline

\multicolumn{1}{|l|}{}          & 25                                   & \multicolumn{1}{c|}{0.4443}   & \multicolumn{1}{c|}{0.2336}    & 0.4783                        & \multicolumn{1}{c|}{0.1712}    & \multicolumn{1}{c|}{0.1064}                  &0.1585                           \\
\multicolumn{1}{|c|}{14}         & 50                               & \multicolumn{1}{c|}{0.2407}    & \multicolumn{1}{c|}{0.1004}     &0.3307                           & \multicolumn{1}{c|}{0.0937}    & \multicolumn{1}{c|}{0.0588}                  &0.1737                           \\
\multicolumn{1}{|l|}{}          & 75                                   & \multicolumn{1}{c|}{0.1411}    & \multicolumn{1}{c|}{0.0321}     &0.2230                           & \multicolumn{1}{c|}{0.0541}    & \multicolumn{1}{c|}{0.0235}                  &0.1821                           \\ \hline

\multicolumn{1}{|l|}{}          & 25                                   & \multicolumn{1}{c|}{0.5529}   & \multicolumn{1}{c|}{0.3160}    & 0.6150                         & \multicolumn{1}{c|}{0.1305}    & \multicolumn{1}{c|}{0.0818}                  & 0.1656                          \\
\multicolumn{1}{|c|}{16}         & 50                               & \multicolumn{1}{c|}{0.3487}    & \multicolumn{1}{c|}{0.1448}     &0.4424                     & \multicolumn{1}{c|}{0.0849}    & \multicolumn{1}{c|}{0.0389}                  & 0.1942                          \\
\multicolumn{1}{|l|}{}          & 75                                  & \multicolumn{1}{c|}{0.1946}    & \multicolumn{1}{c|}{0.0616}     &0.3678                           & \multicolumn{1}{c|}{0.0905}    & \multicolumn{1}{c|}{0.0289}                  & 0.2354                          \\ \hline

\end{tabular}
\end{center}
\caption{Mean and standard deviation of the decreasing of the norm of the residual for the 32 random $L^2(\Omega)$ functions. The results shown are for different number of terms on each of the dimensions tested. }
\label{tab:tableL2}
\end{table}

\subsection{Results for functions in $H^1(\Omega)$} 
\begin{figure}[h!]
\includegraphics[width=7cm]{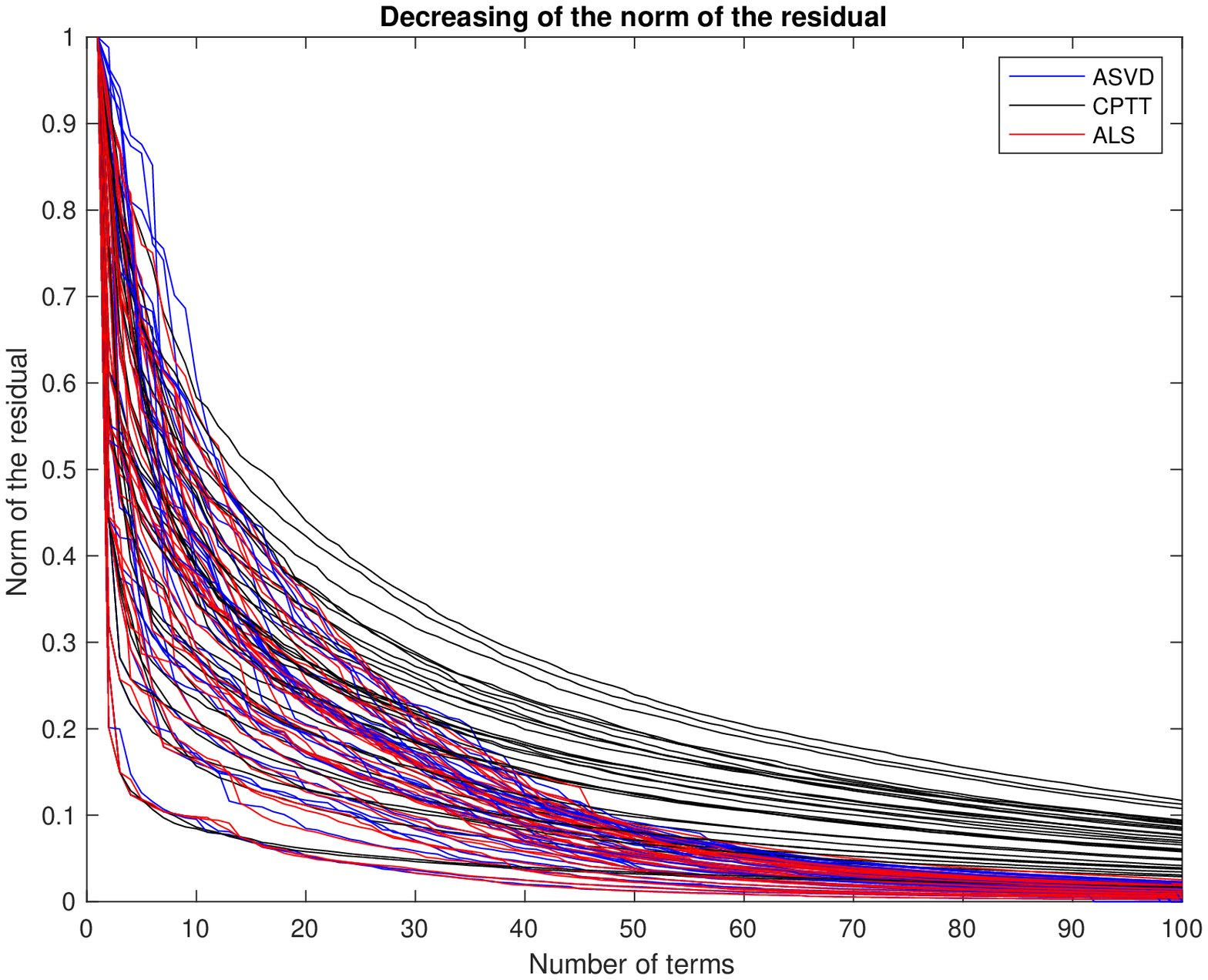}
\includegraphics[width=7cm]{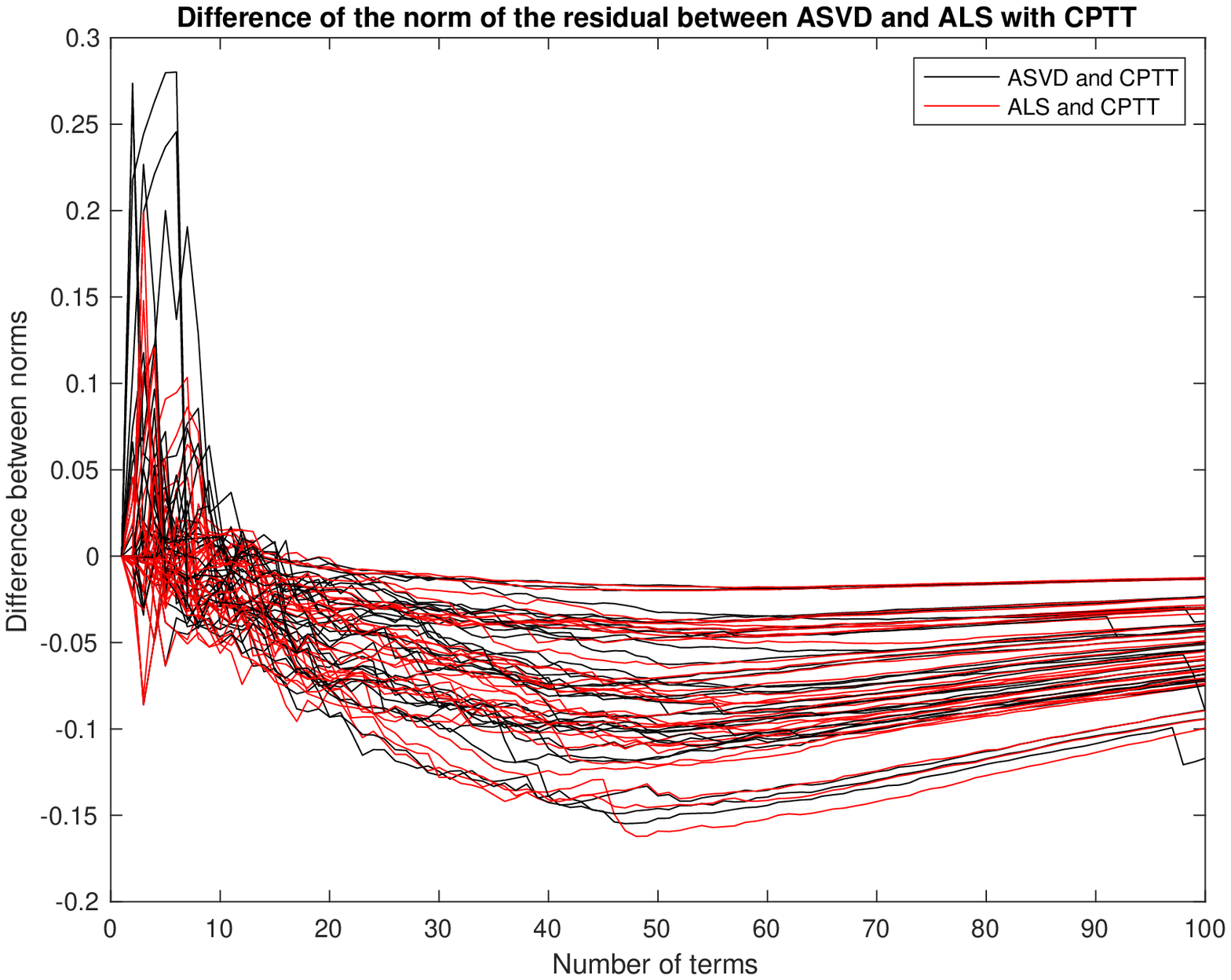}
\caption{Plot associated to functions described in Eq.\eqref{refF}. On the left, the decreasing of the $L^2$ norm of the difference between the exact function $F$ and its approximation in the 4-dimensional tensor of the 32 $H^1$ functions, $\beta=3.1$, using ALS in red, the ASVD method in blue and the CP-TT method in black. On the right, the difference between both norms of the difference between the exact function $F$ and its approximation $\|R\|_{ALS}-\|R\|_{CP-TT}$ in red and $\|R\|_{ASVD}-\|R\|_{CP-TT}$ in black.}
\label{refI7}
\end{figure}

\begin{figure}[h!]
\includegraphics[width=7cm]{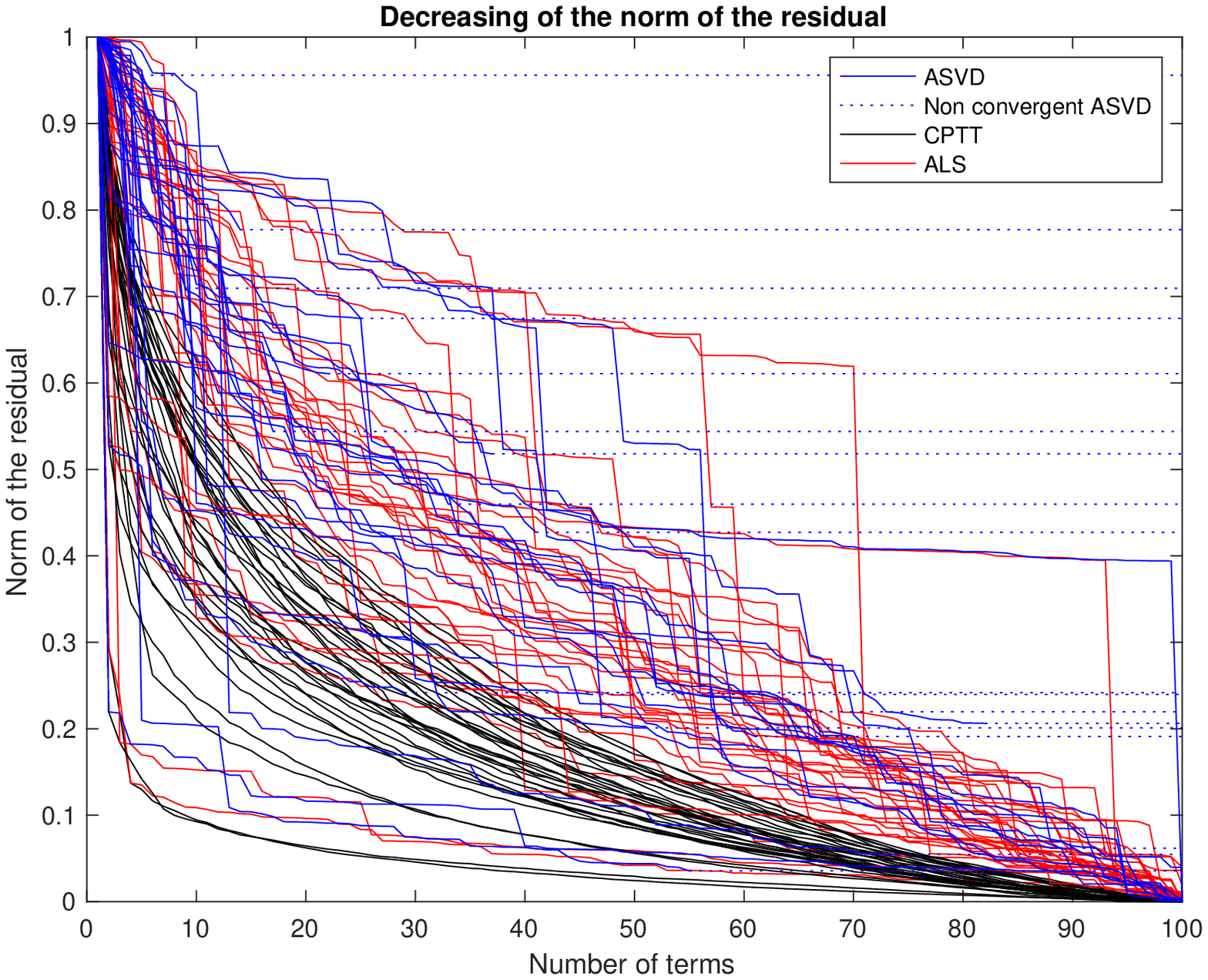}
\includegraphics[width=7cm]{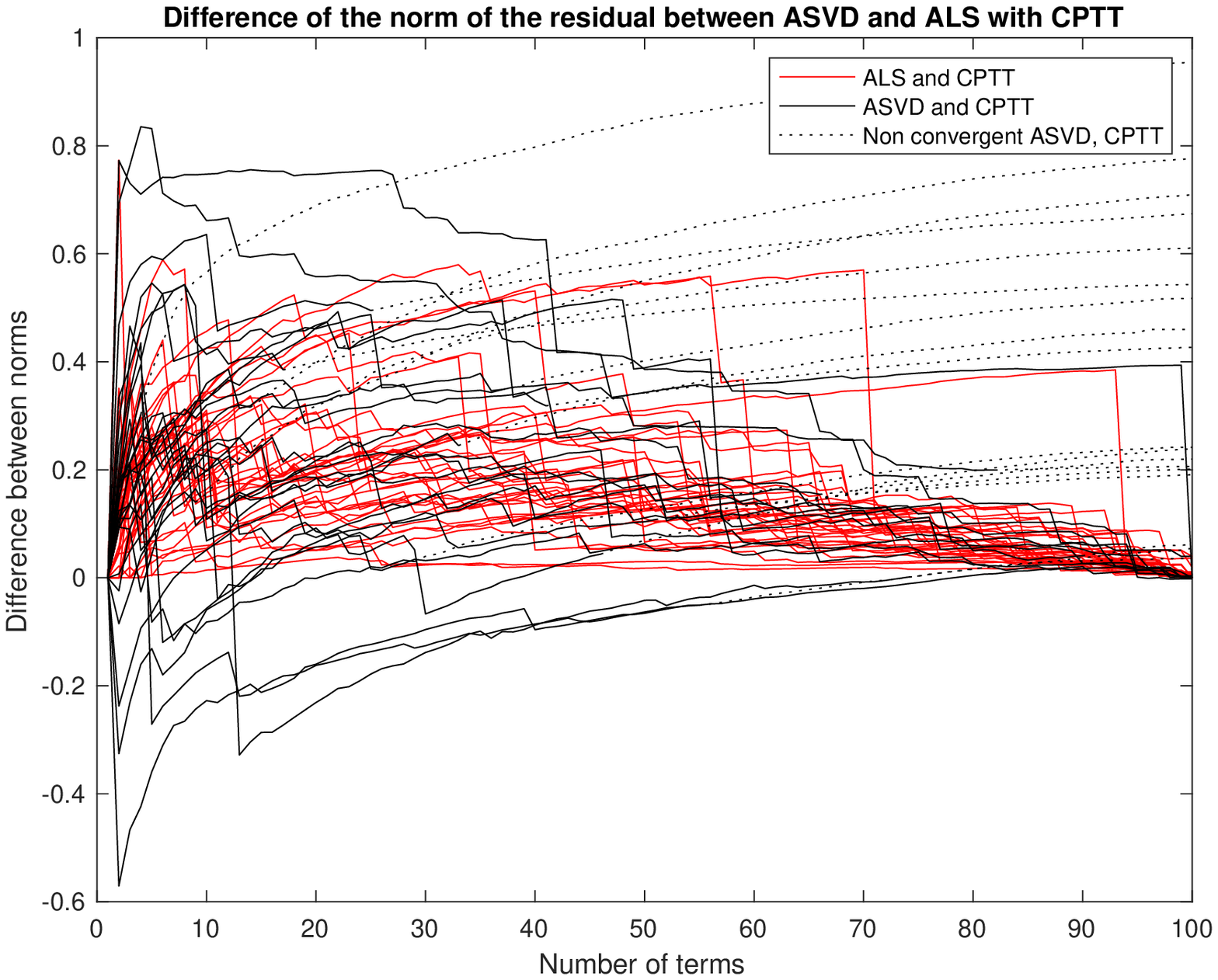}
\caption{Plot associated to functions described in Eq.\eqref{refF}. On the left, the decreasing of the $L^2$ norm of the difference between the exact function $F$ and its approximation in the 12-dimensional tensor of the 32 $H^1$ functions, $\beta=7.1$, using ALS in red, the ASVD method in blue and the CP-TT method in black. On the right, the difference between both norms of the difference between the exact function $F$ and its approximation $\|R\|_{ALS}-\|R\|_{CP-TT}$ in red and $\|R\|_{ASVD}-\|R\|_{CP-TT}$ in black. }
\label{refI8}
\end{figure}

\begin{figure}[h!]
\includegraphics[width=7cm]{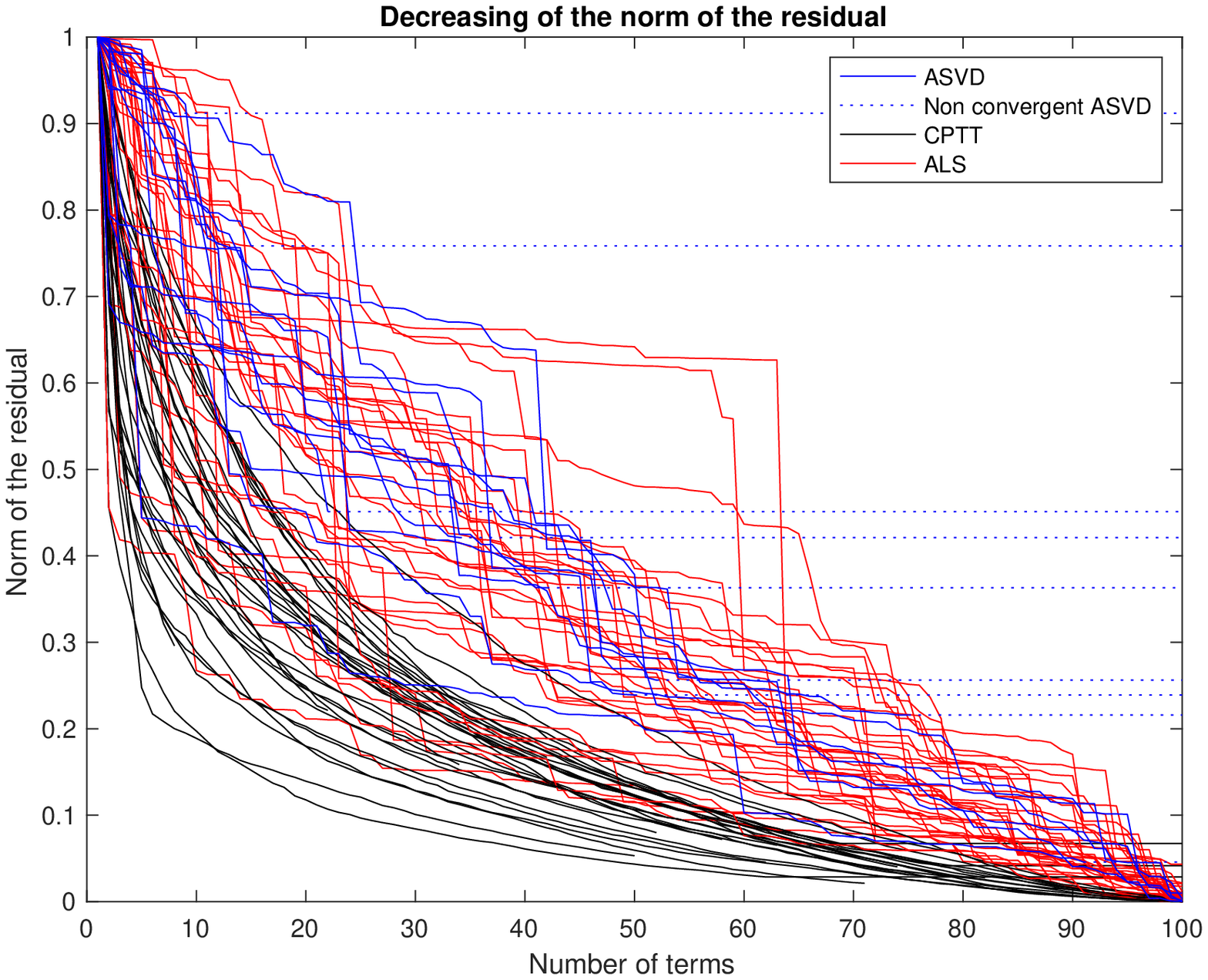}
\includegraphics[width=7cm]{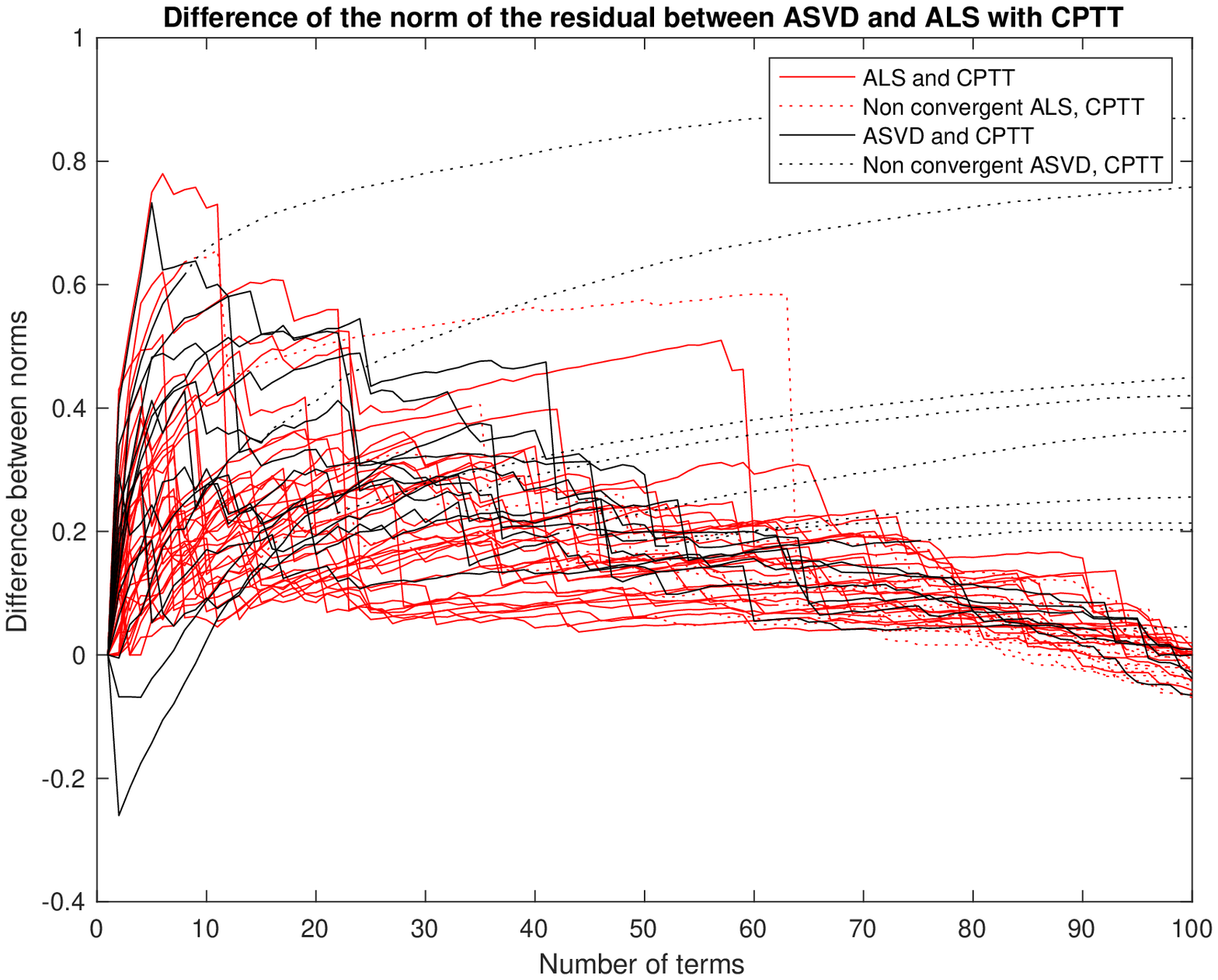}
\caption{Plot associated to functions described in Eq.\eqref{refF}. On the left, the decreasing of the $L^2$ norm of the difference between the exact function $F$ and its approximation in the 16-dimensional tensor of the 32 $H^1$ functions, $\beta=9.1$ using ALS in red, the ASVD method in blue and the CP-TT method in black. On the right, the difference between both norms of the difference between the exact function $F$ and its approximation $\|R\|_{ALS}-\|R\|_{CP-TT}$ in red and $\|R\|_{ASVD}-\|R\|_{CP-TT}$ in black. Both tests are made for 100 terms.}
\label{refI9}
\end{figure}

The method is applied to test cases in which the functions belong to $H^1(\Omega)$, namely when the value of the parameter $\beta$ is chosen to be equal to $\frac{d}{2}+1.1$. The Figures \ref{refI7},\ref{refI8},\ref{refI9}  are the counterpart of Figures  \ref{refI2},\ref{refI3},\ref{refI4} introduced and commented in the previous section. 

The results obtained on $H^1(\Omega)$ functions are equivalent to the ones shown for $L^2(\Omega)$ functions, showing that the decrease in the error norm with the approximation rank is quite regular in CP-TT and behaves in a quite stable way also for higher order tensors.

The \ref{tab:tableH1} collects the results for functions in $H^1(\Omega)$. 
\begin{table}[tbhp]
\begin{center}
\begin{tabular}{cc | ccc | ccc | }
\cline{3-8}
                                &                                      & \multicolumn{1}{c}{}     & \multicolumn{1}{c}{Mean}  & \multicolumn{1}{c|}{}     & \multicolumn{1}{c}{}     & \multicolumn{1}{c}{Standard deviation} & \multicolumn{1}{c|}{}     \\ \hline
\multicolumn{1}{|c|}{Dimension ($d$)} & \multicolumn{1}{c|}{Rank ($r$)} & \multicolumn{1}{c|}{ALS} & \multicolumn{1}{c|}{CPTT} & \multicolumn{1}{c|}{ASVD} & \multicolumn{1}{c|}{ALS} & \multicolumn{1}{c|}{CPTT}              & \multicolumn{1}{c|}{ASVD} \\ \hline
\multicolumn{1}{|l|}{}          & 25                                   & \multicolumn{1}{c|}{0.1722}   & \multicolumn{1}{c|}{0.2261}    &0.1759                          & \multicolumn{1}{c|}{0.0643}    & \multicolumn{1}{c|}{0.2261}                  & 0.1759                       \\
\multicolumn{1}{|c|}{4}         & 50                                 & \multicolumn{1}{c|}{0.0572}   & \multicolumn{1}{c|}{0.1382}     &0.0590                         & \multicolumn{1}{c|}{0.0220}    & \multicolumn{1}{c|}{0.1382}                  & 0.0232                         \\
\multicolumn{1}{|l|}{}          & 75                                   & \multicolumn{1}{c|}{0.0252}    & \multicolumn{1}{c|}{0.0948}    & 0.0262                        & \multicolumn{1}{c|}{0.0103}    & \multicolumn{1}{c|}{0.0948}                  & 0.0110                          \\ \hline

\multicolumn{1}{|l|}{}          & 25                                   & \multicolumn{1}{c|}{0.3741}   & \multicolumn{1}{c|}{0.2938}    &0.4171                         & \multicolumn{1}{c|}{0.1158}    & \multicolumn{1}{c|}{0.0942}                  &0.1341                           \\
\multicolumn{1}{|c|}{6}         & 50                                 & \multicolumn{1}{c|}{0.2037}   & \multicolumn{1}{c|}{0.1507}    &0.2281                         & \multicolumn{1}{c|}{0.0655}    & \multicolumn{1}{c|}{0.0523}                  &0.0791                           \\
\multicolumn{1}{|l|}{}          & 75                                  & \multicolumn{1}{c|}{0.0851}    & \multicolumn{1}{c|}{0.0579}    &0.1045                         & \multicolumn{1}{c|}{0.0334}    & \multicolumn{1}{c|}{0.0233}                  &0.0493                           \\ \hline

\multicolumn{1}{|l|}{}          & 25                                   & \multicolumn{1}{c|}{0.3676}   & \multicolumn{1}{c|}{0.2560}    &0.3977                         & \multicolumn{1}{c|}{0.1361}    & \multicolumn{1}{c|}{0.0905}                  &0.1517                           \\
\multicolumn{1}{|c|}{8}         & 50                                 & \multicolumn{1}{c|}{0.2136}    & \multicolumn{1}{c|}{0.1229}     &0.2413                       & \multicolumn{1}{c|}{0.0807}    & \multicolumn{1}{c|}{0.0451}                  &0.1023                           \\
\multicolumn{1}{|l|}{}          & 75                                   & \multicolumn{1}{c|}{0.1046}    & \multicolumn{1}{c|}{0.0455}     &0.1145                        & \multicolumn{1}{c|}{0.0437}    & \multicolumn{1}{c|}{0.0195}                  &0.0631                           \\ \hline

\multicolumn{1}{|l|}{}          & 25                                  & \multicolumn{1}{c|}{0.4574}   & \multicolumn{1}{c|}{0.3737}    &0.4753                     & \multicolumn{1}{c|}{0.1235}    & \multicolumn{1}{c|}{0.1548}                  & 0.1817                          \\
\multicolumn{1}{|c|}{10}         & 50                                & \multicolumn{1}{c|}{0.2613}    & \multicolumn{1}{c|}{0.3483}     &0.3193                 & \multicolumn{1}{c|}{0.0809}    & \multicolumn{1}{c|}{0.1825}                  &  0.1648                         \\
\multicolumn{1}{|l|}{}          & 75                                   & \multicolumn{1}{c|}{0.1168}    & \multicolumn{1}{c|}{0.3332}     &0.2352                  & \multicolumn{1}{c|}{0.0628}    & \multicolumn{1}{c|}{0.2034}                  & 0.1865                         \\ \hline

\multicolumn{1}{|l|}{}          & 25                                   & \multicolumn{1}{c|}{0.4634}   & \multicolumn{1}{c|}{0.2505}    &0.5182                          & \multicolumn{1}{c|}{0.1681}    & \multicolumn{1}{c|}{0.0842}                  &0.2116                           \\
\multicolumn{1}{|c|}{12}         & 50                               & \multicolumn{1}{c|}{0.2889}    & \multicolumn{1}{c|}{0.1141}     &0.3922                           & \multicolumn{1}{c|}{0.1421}    & \multicolumn{1}{c|}{0.0384}               & 0.2170                          \\
\multicolumn{1}{|l|}{}          & 75                                   & \multicolumn{1}{c|}{0.1278}    & \multicolumn{1}{c|}{0.0382}     &0.3144                           & \multicolumn{1}{c|}{0.0671}    & \multicolumn{1}{c|}{0.0126}               & 0.2502                          \\ \hline

\multicolumn{1}{|l|}{}          & 25                                   & \multicolumn{1}{c|}{0.5943}   & \multicolumn{1}{c|}{0.2169}    &0.4386                         & \multicolumn{1}{c|}{0.2043}    & \multicolumn{1}{c|}{0.1262}                  &0.2014                          \\
\multicolumn{1}{|c|}{14}         & 50                               & \multicolumn{1}{c|}{0.2841}    & \multicolumn{1}{c|}{0.0779}     &0.3132                           & \multicolumn{1}{c|}{0.1277}    & \multicolumn{1}{c|}{0.0686}               &0.1915                           \\
\multicolumn{1}{|l|}{}          & 75                                   & \multicolumn{1}{c|}{0.1422}    & \multicolumn{1}{c|}{0.0244}     &0.2021                           & \multicolumn{1}{c|}{0.0814}    & \multicolumn{1}{c|}{0.0227}               &0.2192                           \\ \hline

\multicolumn{1}{|l|}{}          & 25                                   & \multicolumn{1}{c|}{0.4598}   & \multicolumn{1}{c|}{0.2460}    &0.5543                          & \multicolumn{1}{c|}{0.1496}    & \multicolumn{1}{c|}{0.0726}                  &0.1603                           \\
\multicolumn{1}{|c|}{16}         & 50                               & \multicolumn{1}{c|}{0.2861}    & \multicolumn{1}{c|}{0.1108}     &0.3936                     & \multicolumn{1}{c|}{0.1268}    & \multicolumn{1}{c|}{0.0348}                  &0.2022                           \\
\multicolumn{1}{|l|}{}          & 75                                  & \multicolumn{1}{c|}{0.1395}    & \multicolumn{1}{c|}{0.0438}     &0.3181                        & \multicolumn{1}{c|}{0.0552}    & \multicolumn{1}{c|}{0.0153}                  &0.2477                           \\ \hline

\end{tabular}
\end{center}
\caption{Mean and standard deviation of the decreasing of the norm of the residual for the 32 random $H^1(\Omega)$ functions. The results shown are for different number of terms on each of the dimensions tested. }
\label{tab:tableH1}
\end{table}

\medskip

\paragraph*{Remark:} Despite the fact that one iteration of CP-TT is in general more costly from a computational standpoint with respect to an ALS iteration, at constant target error we computational cost for CP-TT is more or less equivalent (it requires a smaller rank for higher order tensors and the ALS fix point needs more time to converge). Overall, the numerical experiments results are encouraging, showing that CP-TT is a valuable alternative for high-order tensor approximation.


\section{Conclusions and perspectives}
In the present work, a method has been proposed to compute, given a tensor, its CP 
approximation. By leveraging the properties of the TT-SVD algorithm, it is possible to compute a CP decomposition in a stable way, also in the case in which we look for a generic rank$-k$ update.  Albeit the fact that the iterations heavily rely on the TT-SVD algorithm, the proposed strategy does not require to fix \emph{a priori} the order of the variables, but it determines it through an optimization step. 

Several numerical experiments are proposed in order to assess the properties of the method and compare it with ALS and ASVD methods (which are based on similar principles). The experiments suggest that the proposed method can be a valuable tool to compute a CP decomposition in high-dimensional settings, for which the method better behaved with respect to the above mentioned alternatives.  The difference in terms of sparsity between CP-TT and ALS and ASVD methods is large, and even larger for more regular functions, presenting all the CP-TT functions much similar behavior. In order to predict the results without computing them, this has to be taken into account.

The main perspectives of the present work consists in using CP-TT in the solution of multi-linear problems and to investigate how some of the defining steps of this method could be exploited in other tensor formats.

\section*{Acknowledgments}
Virginie Ehrlacher acknowledges support from the ANR COMODO project (ANR-19-CE46-0002).

Damiano Lombardi acknowledges support from the ANR ADAPT project (ANR-18-CE46-0001).

This publication is part of a project that has received funding from the European Research Council (ERC) under the European Union’s Horizon 2020 Research and Innovation Programme – Grant Agreement $n^{\circ}$ 810367.

\bibliographystyle{plain}
\bibliography{biblioCPTT}

\end{document}